\batchmode
\documentclass[11pt]{article}

\usepackage{epsfig}
\usepackage{graphicx}
\usepackage{color}
\usepackage{mathtools}

\newtheorem{theorem}{Theorem}
\newtheorem{lemma}{Lemma}
\newtheorem{corollary}{Corollary}

\newcommand{\be}{\begin{equation}}
\newcommand{\ee}{\end{equation}}
\newcommand{\bea}{\begin{eqnarray}}
\newcommand{\eea}{\end{eqnarray}}
\newcommand{\beas}{\begin{eqnarray*}}
\newcommand{\eeas}{\end{eqnarray*}}
\newcommand{\ba}{\begin{array}}
\newcommand{\ea}{\end{array}}

\DeclarePairedDelimiter{\floor}{\lfloor}{\rfloor}
\DeclarePairedDelimiter{\ceil}{\lceil}{\rceil}

\definecolor{armygreen}{rgb}{0.29, 0.33, 0.13}

\newcommand{\real}{\mbox{$\mathbb{R}$}}

\def\XXint#1#2#3{{\setbox0=\hbox{$#1{#2#3}{\int}$}
     \vcenter{\hbox{$#2#3$}}\kern-.5\wd0}}

\newcommand{\mcL}{\ensuremath{\mathcal{L}}}
\newcommand{\mcM}{\ensuremath{\mathcal{M}}}
\newcommand{\mcN}{\ensuremath{\mathcal{N}}}

\newcommand{\mcS}{\ensuremath{\mathcal{S}}}

\newcommand{\mrI}{\ensuremath{\mathrm{I}}}

\newcommand{\wtilde}{\ensuremath{\widetilde}}
\newcommand{\what}{\ensuremath{\widehat}}

\def\qed{\hbox{\vrule width 6pt height 6pt depth 0pt}}

\usepackage{amsmath}
\usepackage{amssymb}

\title{Analysis and Petrov-Galerkin  numerical approximation for variable coefficient two-sided fractional diffusion, 
advection, reaction equations} 
\author{
	Xiangcheng Zheng \thanks{School of Mathematical Sciences, 
	Peking University, Beijing 100871, China. 
	email: {\tt zhengxch@math.pku.edu.cn}.} 
	\and
	V.J.~Ervin\thanks{School of Mathematical and Statistical Sciences,
	  Clemson University, Clemson, South Carolina 29634-0975, USA.
	  email: {\tt vjervin@clemson.edu}. }
	\and 
	 Hong Wang \thanks{Department of Mathematics, University of South Carolina, Columbia,
		South Carolina 29208, USA. email: {\tt  hwang@math.sc.edu}.} 
 }

\date{\today}

\begin{document}
\maketitle

\begin{abstract}
In this paper we investigate the variable coefficient two-sided fractional diffusion, advection, reaction equations on a bounded interval. It is known that the fractional diffusion operator may lose coercivity due to the variable coefficient, which makes both the mathematical and numerical analysis challenging. To resolve this issue, we design appropriate test and trial functions to prove the inf-sup condition of the variable coefficient fractional diffusion, advection, reaction operators in suitable function spaces. Based on this property, we prove the well-posedness and regularity of the solutions, as well as analyze the Petrov-Galerkin approximation scheme for the proposed model.
Numerical experiments are presented to substantiate the theoretical findings and to compare the behaviors of different models.
\end{abstract}



\setcounter{equation}{0}
\setcounter{figure}{0}
\setcounter{table}{0}
\setcounter{theorem}{0}
\setcounter{lemma}{0}
\setcounter{corollary}{0}
\setcounter{definition}{0}
\section{Introduction}
 \label{sec_intro}
In this work we investigate the well-posedness, regularity and numerical approximation of the variable coefficient two-sided fractional diffusion, advection, reaction equation
\begin{align}
  \wtilde{\mcL}_{r}^{\alpha}u(x) \ + \ b(x) D u(x) \ + \ c(x) u(x)   &= \ f(x) \, , \ \  x \in \mrI \, , 
   \label{DefProb2}  \\
   \mbox{subject to } u(0) \, = \, u(1) &= \,  0.
\label{DefBC2} 
\end{align}  
Here $\mrI \, := \, (0 , 1)$,  $b \in W^{1,\infty}(\mrI)$, $c,k \in L^{\infty}(\mrI)$ with $k(x)\geq k_0>0$ represent advection, reaction, and material diffusivity coefficients, respectively, $D$ denotes the usual differential operator,  The fractional diffusion operator $\wtilde{\mcL}_{r}^{\alpha}$ with $1 < \alpha < 2$ and $0 \le r \le 1$ has two possible forms $\acute{\mcL}_{r}^{\alpha}$ and $\grave{\mcL}_{r}^{\alpha}$ according to the position of the variable diffusivity coefficient 
\begin{align}
\acute{\mcL}_{r}^{\alpha}u(x) &:= \  
   - D  \big(\big( r D^{-(2 - \alpha)} \ + \ (1 - r) D^{-(2 - \alpha)*}  \big) k(x) \, D u(x) \big) \, ,   \label{defmcL}  \\
\grave{\mcL}_{r}^{\alpha}u(x) &:= \ 
   - D \big(  k(x) \big( r D^{-(2 - \alpha)} \ + \ (1 - r) D^{-(2 - \alpha)*}  \big) \, D u(x)  \big) ,   \label{defmcL2}  
\end{align}
where $D^{-(2-\alpha)}$ and $D^{-(2-\alpha)*}$ refer to the $(2-\alpha)$-order left and right fractional integral operators, respectively, defined by
\begin{align}
D^{-(2-\alpha)} u(x) &:=  \frac{1}{\Gamma(2 - \alpha)}   \int_{0}^{x} \frac{1}{(x - s)^{\alpha - 1}} \, u(s) \, ds \, ,   \label{defDalpha}  \\
D^{-(2-\alpha)*} u(x) &:=   \frac{1}{\Gamma(2 - \alpha)} \int_{x}^{1} \frac{1}{(s - x)^{\alpha - 1}} \,  u(s) \, ds \, .   \label{defDalpha*} 
\end{align}

There are fruitful mathematical and numerical  results for space-fractional differential equations in the literature \cite{che161,che162,gin191,jin151,LiTel,liu041,mao182,She20SINUM,wang21jcam,xu141,zay131}. The well-posedness of a constant diffusivity coefficient analogue of problem \eqref{DefProb2}-\eqref{DefBC2} and error estimates of its Galerkin finite element approximation was proved in \cite{erv061}. It was shown in \cite{wan131} that Galerkin weak formulation of the diffusive only version (i.e., $b=c=0$)  of problem \eqref{DefProb2}-\eqref{DefBC2} is not coercive and the corresponding Galerkin finite element approximation is not guaranteed to converge \cite{wan152}. It was also shown in \cite{jin151,wan152a} that the solution to the diffusive only version of problem \eqref{DefProb2}-\eqref{DefBC2} or even its constant diffusivity coefficient analogue exhibits weak singularity near the end points of the interval, which makes the full regularity assumption on the true solution of fractional differential equations in many optimal-order error estimates of numerical approximations in the literature unrealistic. An optimal-order error estimate of a spectral Galerkin method in the $L^2$ norm was proved in \cite{wan151} for a one-sided diffusive only version of problem \eqref{DefProb2}-\eqref{DefBC2} by assuming only the smoothness of the diffusivity coefficient and the right-hand side function. 
The complete regularity analysis of the diffusive only version of problem \eqref{DefProb2}-\eqref{DefBC2} was given in \cite{erv162}, which were subsequently extended to the constant diffusivity coefficient analogue of problem \eqref{DefProb2}-\eqref{DefBC2} in \cite{erv191} (see also \cite{hao201}).

Despite the aforementioned progress, the corresponding results for the variable coefficient fractional diffusion, advection, reaction equations are largely missing. There exists some recent work on Petrov-Galerkin approximations to two-sided fractional diffusion, reaction equations \cite{hao202,hao21anm}, one-sided fractional diffusion, advection, reaction equations \cite{erv17,fud21,jin162}, and two-sided fractional diffusion, advection, reaction equations \cite{zhe211}, all with constant diffusivity coefficients. To the best of our knowledge, the only available result for the Petrov-Galerkin method applied to variable coefficient fractional diffusion problems is \cite{wan131}, in which the weak coercivity in the sense of inf-sup condition was proved for the one-sided variable coefficient fractional diffusion operator, i.e.,  $\grave{\mcL}_{r}^{\alpha}$ with $r=1$. Then a discontinuous Petrov-Galerkin method was applied to this model in \cite{wan152} for numerical approximation. How to prove the weak coercivity for the two-sided case (\ref{DefProb2}) remains untreated in the literature, which hinders the theoretical and numerical analysis. 

In this paper we analyze the well-posedness, regularity and numerical approximation of problem (\ref{DefProb2})-(\ref{DefBC2}) by designing appropriate test and trial functions to prove the weak coercivity of the Petrov-Galerkin weak formulation for problem (\ref{DefProb2})-(\ref{DefBC2}). For different choices of $\wtilde{\mcL}_{r}^{\alpha}$ we select different test and trial spaces for the sake of the proof. Particular attention has been paid in the choice of the weak formulations given, in order that the 
required regularity for $k(x)$, $b(x)$ and $c(x)$ are consistent with those for the usual (second-order) diffusion, advection, reaction problem. Error estimates of the Petrov-Galerkin spectral approximation scheme are proved in both weighted $L^2$ and energy norms, which provides theoretical supports for numerical computations.

The rest of the paper is organized as follows. In Section 2 we introduce definitions, notations, and useful results. In Section 3 we analyze properties of the Petrov-Galerkin weak formulation for
\eqref{DefProb2}-\eqref{DefBC2}, and establishes the existence and uniqueness of its solution. In Section 4 we prove the regularity of the solutions.
The Petrov-Galerkin approximation scheme is proposed and analyzed in Section 5.
Numerical experiments are presented in the last section to substantiate the theoretical findings and to compare the behaviors of different models.

 \setcounter{equation}{0}
\setcounter{figure}{0}
\setcounter{table}{0}
\setcounter{theorem}{0}
\setcounter{lemma}{0}
\setcounter{corollary}{0}
\setcounter{definition}{0}
\section{Notation and properties}
\label{sec_not}
In this section, we present various notation and spaces to be used subsequently. We let $\mathbb{N}_{0}  := \mathbb{N} \cup \{0\}$ and  use $y_{n} \sim n^{p}$ to denote that there exist constants $c$ and $C > 0$ such that, as 
 $n \rightarrow \infty$,  
 $c \, n^{p} \le | y_{n} | \le C \, n^{p}$. Additionally, we use $a \, \lesssim \, b$ and $a \, \simeq \, b$
 to denote that there exists constants $C_{1}$ and $C_{2}$ such that
 $a \, \le \, C_{2} \,  b$, and $C_{1} b \, \le \, a \, \le \, C_{2} \, b$, respectively. 
 For $t \in \mathbb{R}$, $\floor{t}$ is used to denote the largest integer that is less than or equal to $t$, and
 $\ceil{t}$ is used to denote the smallest integer that is greater than or equal to $t$.
\subsection{Jacobi polynomials}
Jacobi polynomials have close relations to the fractional problems
\cite{aco181, erv162, mao181, mao161}. The classical Jacobi (orthogonal) polynomials $\{P_{n}^{(a , b)}(x)\}_{n\geq 0}$ for $n \ge 0$ and $a, b > -1$ are defined on $[-1,1]$ (see \cite{abr641,sze751}).
%
As we are interested in the domain $\mrI = (0 , 1)$, we let $G_{n}^{(a , b)}(x) \, = \, P_{n}^{(a , b)}( 2 x \, - \, 1 )$.
Then
\begin{align}
\int_{0}^{1}   (1 - x)^{a}  \, x^{b} \, G_{j}^{(a , b)}(x) \, G_{k}^{(a , b)}(x)  \,  dx 
 &= \
   \left\{ \begin{array}{ll} 
   0 , & k \ne j \, , \\
   |\| G_{j}^{(a , b)} |\|^{2}
   \, , & k = j  
   \, . \end{array} \right.    \nonumber \\
 \quad \quad \mbox{where } \  \ |\| G_{j}^{(a , b)} |\| \ = \ |\| G_{j}^{(b , a)} |\| &= \
 \left( \frac{1}{(2j \, + \, a \, + \, b \, + 1)} 
   \frac{\Gamma(j + a + 1) \, \Gamma(j + b + 1)}{\Gamma(j + 1) \, \Gamma(j + a + b + 1)}
   \right)^{1/2} \, .  \label{spm22g} 
\end{align}        

From \cite[equation (2.9)]{zhe20anm} we have
\begin{equation}\label{GnNormRelation}
\frac{1}{2} \, \le \, \frac{  |\| G_{j}^{(a - b \, , \, b)}  |\|^{2} }{  |\| G_{j + 1}^{(b - 1 \, , \, a - b - 1)} |\|^{2}} 
	\, = \, \frac{j + 1}{j + a} \, \le \, 1, \, \quad j \ge 0.
\end{equation}

 From \cite{mao161} it follows that
\begin{align}
\frac{d^{k}}{dt^{k}} G_{n}^{(a , b)}(t) 
  &= \ \frac{\Gamma(n + k + a + b + 1)}{  \Gamma(n + a + b + 1)} 
  G_{n - k}^{(a + k \, , \, b + k)}(t)  \, .
    \label{eqC4}  
 \end{align}

Note that, from Stirling's formula, we have that
\begin{equation}
\lim_{n \rightarrow \infty} \, \frac{\Gamma(n + \sigma)}{\Gamma(n) \, n^{\sigma}}
\ = \ 1 \, , \mbox{ for } \sigma \in \mathbb{R}.  
 \label{eqStrf}
\end{equation} 
 
For compactness of notation, let
 $\omega^{(a , b)} \, = \, \omega^{(a , b)}(x) \, := \, (1 - x)^{a} \, x^{b} $.

A useful formula used in the analysis below is \cite{mao161}
\begin{align}
\frac{d^{k}}{dx^{k}} \omega^{(a + k \, , \, b + k)}(x) \, G_{n - k}^{(a + k \, , \, b + k)}(x) 
  &= \ \frac{(-1)^{k} \, n !}{ (n - k) !} \, 
  \omega^{(a  ,  b)}(x) \, G_{n}^{(a , b)}(x)   \, , \mbox{ for } 0 \le k \le n \, .  \label{eqC4p5}
\end{align}   

\textbf{Condition A}:
The parameters $\alpha$, $\beta$, $r$ and constant $c_{*}^{*}$ satisfy:
$1 < \alpha < 2$, $\alpha - 1 \, \le  \, \beta \, , \, \alpha - \beta \, \le  \, 1$, $0 \le r \le 1$ 
\begin{equation}
  c_{*}^{*} \ = \ \frac{\sin(\pi \alpha)}{\sin(\pi (\alpha - \beta)) \, + \, \sin(\pi \beta)} \, < 0 \, ,  \label{defcss}
\end{equation}  
where $\beta$ is determined by
\begin{equation}
  r \ = \ \frac{\sin( \pi \, \beta)}{\sin( \pi ( \alpha - \beta)) \, + \,  \sin( \pi \, \beta)} \, . \label{propker0} 
\end{equation}

From \cite{jia181}, we have for $\mcM^{*} v(x)\, := \, \big( r D^{-(2 - \alpha)*} \ + \ (1 - r) D^{-(2 - \alpha)}  \big) D v(x)$ 
\begin{align}
   \mcM^{*}  \omega^{(\beta , \, \alpha - \beta)}(x) \, G_{k}^{(\beta \, , \, \alpha - \beta)}(x)
   &= \ \mu_{k} \, G_{k+1}^{(\alpha - \beta - 1 \, , \, \beta - 1)}(x) \, ,  \label{propmcM}  \\
   \mbox{ where } \ \ 
   \mu_{k} &= \ c_{*}^{*}
    \frac{\Gamma(k + \alpha)}{\Gamma(k + 1)}  \, ,
    \ k = 0, 1, 2, \ldots,
\label{propmcM1}  \\
 \mbox{Also, using \eqref{eqStrf} },  \  | \mu_{k} | &\sim \, (k + 1)^{\alpha - 1} \mbox{  as  } \ k \rightarrow \infty . 
  \label{propmcM2} 
\end{align}
 
Using \eqref{eqC4p5} we have
\be
D \omega^{(\alpha - \beta , \, \beta)}(x) \,   G_{n}^{(\alpha - \beta  ,  \beta)}(x) \ = 
\ - (n + 1) \,  \omega^{(\alpha - \beta - 1 \, , \, \beta - 1)}  \, G_{n+1}^{(\alpha - \beta - 1 \, , \, \beta - 1)}(x) \, , 
\ n = 0, 1, 2, \ldots
 \label{propMnx}
\ee
\subsection{Function spaces}
For $\sigma(x) > 0, \ x \in (0 , 1)$, let 
\begin{equation}
L_{\sigma}^{2}(\mrI) \, := \, \left \{ f(x) \, : \, \int_{0}^{1} \sigma(x) \, f(x)^{2} \, dx \ < \ \infty \right \} \, .
\label{defLw}
\end{equation}
Associated with $L_{\sigma}^{2}(0 , 1)$ is the inner product, $( \cdot , \cdot )_{\sigma}$, and
norm, $\| \cdot \|_{\sigma}$, defined by
\[
( f \,  , \, g )_{\sigma} \, := \, \int_{0}^{1} \sigma(x) \, f(x) \, g(x) \, dx \, , \quad \mbox{and} \quad
 \| f \|_{\sigma} \, := \, \left( \big ( f \,  , \, f \big )_{\sigma} \right)^{1/2} \, .
\]
Without a subscript, $( \cdot , \cdot )$ denotes the usual $L^{2}(\mrI)$ inner product.
The set of orthogonal polynomials $\{ G_{j}^{(a , b)} \}_{j = 0}^{\infty}$ form an orthogonal basis
for $L^{2}_{\omega^{(a , b)}}(\mrI)$, and for $\widehat{G}_{j}^{(a , b)} \, := \, G_{j}^{(a , b)} / |\| G_{j}^{(a , b)} |\|$, 
$\{ \what{G}_{j}^{(a , b)} \}_{j = 0}^{\infty}$ form an orthonormal basis
for $L^{2}_{\omega^{(a , b)}}(\mrI)$.

The weighted Sobolev spaces $H^{s}_{(a , b)}(\mrI)$ could be defined via two different but equivalent ways. In the first
definition, for the non-negetive integer $s$, define $H^s_{\omega^{(a,b)} }(\mrI)$ for $a, b > -1$ as \cite{bab011,guo041}
\begin{equation} 
H^s_{\omega^{(a,b)} }(\mrI) := \bigg \{ v \, : \| v \|_{s , \omega^{(a,b)}}^2 := 
 \sum_{j=0}^s \big \| D^j v \big \|_{\omega^{(a+j,b+j)}}^2 < \infty \bigg \} .
\label{defHw}
\end{equation} 
The spaces $H^s_{\omega^{(a,b)} }(\mrI)$ for $s > 0 \not \in \mathbb{N}$ are defined by the
$K$- method of interpolation, while for $s < 0$
the spaces are defined by (weighted) $L^{2}$ duality. The second definition is based on the decay rate of the Jacobi 
coefficients of a function. For $v \in L^{2}_{\omega^{(a , b)}}(\mrI)$, the following expansion holds
\be
v(x) \ = \ \sum_{j = 0}^{\infty} v_{j} \, \what{G}_{j}^{(a , b)}(x) ,~~  v_{j} \ = \ \int_{0}^{1}  \omega^{(a , b)}(x) \, v(x) \,  \what{G}_{j}^{(a , b)}(x) \, dx \, .
 \label{defvj}
\ee
Then for $s, a, b \in \real$,  $a, b > -1$, 
$L^{2}_{(a , b)}(\mrI) \, := \, L^{2}_{\omega^{(a , b)}}(\mrI)$, define 
\be
H^{s}_{(a , b)}(\mrI) \, := \, \left \{v \, : \,
 \sum_{j = 0}^{\infty} (1 + j^{2})^{s} \, v_{j}^{2} \, < \, \infty \right \}.
\label{defHr}
\ee

\begin{theorem}  \cite[Theorem 4.1]{erv191} \label{thmeq2}
The spaces $H^{s}_{(a , b)}(\mrI)$ and $H^{s}_{\omega^{(a , b)}}(\mrI)$
coincide, and their corresponding norms are equivalent. 
\end{theorem}


\begin{lemma}  \cite[Lemma 4.5]{erv191}\label{lmamapD}
For $s, a, b \in \real$,  $a, b > -1$, the differential operator $D$ is a bounded mapping from $H^{s}_{(a , b)}(\mrI)$ onto
$H^{s - 1}_{(a + 1 \, , \, b + 1)}(\mrI)$.
\end{lemma}

For convenience, from  hereon we use $H^{s}_{(a , b)}(\mrI)$ 
to represent the spaces
$H^s_{\omega^{(a,b)} }(\mrI)$ and $H^{s}_{(a , b)}(\mrI)$.

Let $\mcS_{N}$ denote the space of polynomials of degree  less than or equal to $N$. We define the weighted $L^2$ orthogonal projection
$P_{N} : \, L^{2}_{\omega}(\mrI) \rightarrow \mcS_{N}$ by the condition
\begin{equation}\label{Proj}
\big ( v \, - \,  P_{N}v \ , \ \phi_N \big)_{\omega} \ = \ 0 \, , \ \ \forall \phi_N \in \mcS_{N}.
\end{equation}
Note that $P_{N}v \ = \ \sum_{j = 0}^{N} v_{j} \, \what{G}_{j}^{(a , b)}(x)$, 
where $v_{j}  \ = \ \int_{0}^{1}  \omega(x) \, v(x) \,  \what{G}_{j}^{(a , b)}(x) \, dx$.
\begin{lemma}\label{lem:Approx} \cite[Theorem 2.1]{guo041}
For $\mu \in \mathbb{N}_{0}$ and $v \in H^{t}_{\omega}(\mrI)$, with $0 \le \mu \le t$, there exists a
constant $C$, independent of $N, \, \alpha$ and $\beta$ such that
\begin{equation}\label{Approx}
\big \| v  -  P_{N} v \|_{H^{\mu}_{\omega}(\mrI)} \ \le \ C \, 
N^{\mu - t} \, \| v \|_{H^{t}_{\omega}(\mrI)}.
\end{equation}
\end{lemma}  

\textbf{Remark}: In \cite{guo041} \eqref{Approx} is stated for $t \in \mathbb{N}_{0}$. The result extends to 
$t \in \mathbb{R}^{+}$ using interpolation.



Define the space $W^{k , \infty}_{w}(\mrI)$ for $k \in \mathbb{N}_{0}$ and its associated norm as
\begin{align}
W^{k , \infty}_{w}(\mrI) &:= \ \left\{ f \, : \ (1 - x)^{j/2} x^{j/2} D^{j}f(x) \in L^{\infty}(\mrI) , \ j = 0, 1, \ldots, k \right\} , 
\label{defWinftw}  \\
\| f \|_{W^{k , \infty}_{w}(\mrI)} &:= \ \max_{0 \le j \le k} \| (1 - x)^{j/2} x^{j/2} D^{j}f(x) \|_{L^{\infty}(\mrI)} \, .
\label{defWinftwnorm}
\end{align}
The subscript $w$ denotes that $W^{k , \infty}_{w}(\mrI)$ is a weaker space than $W^{k , \infty}(\mrI)$. 
\begin{lemma} \label{lmaprodsp}  \cite[Lemma 7.1]{erv191}
Let $s \ge 0$, $\alpha, \, \beta > -1$, $k \ge s$, and $f \in W^{k , \infty}_{w}(\mrI)$. For
$g \in H^{s}_{(\alpha , \beta)}(\mrI)$
\be
\ f g \in H^{s}_{(\alpha , \beta)}(\mrI) \ \mbox{ and   } 
\| f g \|_{H^{s}_{(\alpha , \beta)}(\mrI)} \ \le \  \| f \|_{W^{k , \infty}_{w}(\mrI)} \, 
 \| g \|_{H^{s}_{(\alpha , \beta)}(\mrI)} \, .
   \label{prodr1}  
\ee
\end{lemma}

For $k-1 < s < k$ the requirement that $f \in W^{k , \infty}_{w}(\mrI)$ can be slightly relaxed to
$f \in W^{k - 1 \,  , \, \infty}_{w}(\mrI)$ \underline{and} 
$| f^{(k -1)}(x) \, - \, f^{(k -1)}(y) | \ \le \ C \, | x \, - \, y |^{\sigma}$, 
for $\sigma > \, s - \, (k - 1)$ (see the proof of Lemma 5.3 in \cite{din121}).

For compactness of notation, for $\alpha$ and $r$ defined in \eqref{DefProb2} and $\beta$ defined in 
\eqref{propker0} we introduce
\begin{equation}
 \omega(x) \, := \, \omega^{(\alpha - \beta , \beta)}(x) \, = \, (1 - x)^{\alpha - \beta} \, x^{\beta} \, ,  \ 
 \mbox{ and } \ \omega^{*}(x) \, := \, \omega^{(\beta  , \alpha - \beta)}(x) \, = \, (1 - x)^{\beta} \, x^{\alpha - \beta} \, .  \label{defomega}
\end{equation}

The following lemma is helpful in establishing the error estimate between $u$ and its approximation $u_{N}$.
(See Corollary \ref{uerrL2}.)
\begin{lemma} \label{lma4Her} \cite[Lemma 4.3]{zhe211}
Let $0 \le \mu \le 1$. For $\zeta \in H^{\mu}_{\omega}(\mrI)$, then $z \, := \, \omega \, \zeta \in 
H^{\mu}_{\omega^{-1}}(\mrI)$, with, for some $C > 0$,
\be
 \| z \|_{H^{\mu}_{\omega^{-1}}(\mrI)} \ \le \ C \,  \| \zeta \|_{H^{\mu}_{\omega}(\mrI)} \, . \label{zestH}
\ee
\end{lemma}

In the rest of the paper we use $\langle \cdot , \cdot \rangle_{\omega}$ to denote the weighted $L^{2}$ duality
pairing between functions if $H^{-s}_{(\alpha - \beta \, , \, \beta)}(\mrI)$ and 
$H^{s}_{(\alpha - \beta \, , \, \beta)}(\mrI)$.

 \setcounter{equation}{0}
\setcounter{figure}{0}
\setcounter{table}{0}
\setcounter{theorem}{0}
\setcounter{lemma}{0}
\setcounter{corollary}{0}
\setcounter{definition}{0}
\section{Analysis of weak formulation}
\label{sec_form}
We mainly consider model \eqref{DefProb2}-\eqref{DefBC2} with $\wtilde{\mcL}_{r}^{\alpha}=\acute{\mcL}_{r}^{\alpha}$, and briefly address the case of $\wtilde{\mcL}_{r}^{\alpha}=\grave{\mcL}_{r}^{\alpha}$. (Since the latter is the adjoint of the former, the proofs are similar.)

\underline{Weak Formulation} \\
 Given $f \in H^{-(\alpha - 1)}_{\omega^{*}}(\mrI)$,  $0 < k_{0} \le k(x)  \in L^{\infty}(\mrI)$,
$b \in W^{1 , \infty}_{w}(\mrI)$ and 
$c \in L^{\infty}(\mrI)$, 
determine $\phi \in H^{1}_{\omega}(\mrI)$
such that $u(x) \ = \ \omega(x) \, \phi(x)$ satisfies
\be
\langle \acute{\mcL}_{r}^{\alpha} u \ + \ b \, D u \ +  \ c \, u \, , \, \psi \rangle_{\omega^{*}} 
\ = \ \langle f \, , \, \psi \rangle_{\omega^{*}} \, , \ \ \forall \, \psi \in H^{\alpha - 1}_{\omega^{*}}(\mrI) \, .
\label{wform1}
\ee
For simplicity of the notations, let $B \, : \, H^{1}_{\omega}(\mrI) \times H^{\alpha - 1}_{\omega^{*}}(\mrI) \rightarrow \real$ and
$F \, :  \, H^{\alpha - 1}_{\omega^{*}}(\mrI) \rightarrow \real$ be defined by
\begin{align}
 B(\phi , \psi) &:= \ \langle \acute{\mcL}_{r}^{\alpha} u \ + \ b \, D u \ +  \ c \, u \, , \, \psi \rangle_{\omega^{*}} \, =: \, 
 B_0(\phi , \psi) \ + \ B_1(\phi , \psi) \ + \ B_2(\phi , \psi),
 \label{defB}  \\
F(\psi) &:= \ \langle f \, , \, \psi \rangle_{\omega^{*}} \, .   \label{defF}
\end{align}

We recall the
Banach-Ne\v{c}as-Babu\v{s}ka theorem to support the subsequent analysis.
\begin{theorem} \cite[Pg. 85, Theorem 2.6]{ern041} \label{BNBthm}
Let $H_{1}$ and $H_{2}$ denote two real Hilbert spaces, $B(\cdot , \cdot) \, : \, H_{1} \times H_{2} \rightarrow \real$
a bilinear form, and $F \, : \, H_{2} \rightarrow \real$ a bounded linear functional on $H_{2}$. Suppose there are constants
$C_{1} < \infty$ and $C_{2} > 0$ such that
\begin{align}
&(i) \ | B(w , v) | \, \le \, C_{1} \, \| w \|_{H_{1}} \, \| v \|_{H_{2}} \, ,  \ 
\mbox{ for all } w \in H_{1} \, , \ v \in H_{2} \, , \label{BB1}  \\
&(ii) \ \sup_{0 \ne v \in H_{2}} \frac{ | B(w , v) |}{ \| v \|_{H_{2}} } \ \ge \ C_{2} \,  \| w \|_{H_{1}} \, ,  
 \mbox{ for all } w \in H_{1} \, , \label{BB2}  \\
&(iii) \ \sup_{w \in H_{1}}  | B(w , v) | \ > \ 0 \, ,  \mbox{ for all } v \in H_{2} \, , \ v \neq 0 \, .  \label{BB3}  
\end{align}
Then there exists a unique solution $w_{0} \in H_{1}$ satisfying $B(w_{0} \, , \, v) \, = \, F(v)$ for all 
$ v \in H_{2}$. Further, $\| w_{0} \|_{H_{1}} \, \le \, \frac{1}{C_{2}} \| F \|_{H_{2}}$.
\end{theorem}


\subsection{Property \eqref{BB1} -- Continuity of $B(\cdot , \cdot)$}
\label{ssec_Bcty}
In this section we establish property \eqref{BB1}.

\underline{Continuity of $B_{0}(\cdot , \cdot)$} \\
To establish the continuity of $B_{0}(\cdot , \cdot)$, we begin by choosing 
$\phi(x) \ = \ \sum_{i = 0}^{\infty} \phi_{i} \, \what{G}_{i}^{(\alpha - \beta , \beta)}(x) \, \in H^{1}_{\omega}(\mrI)$. Then,
\be
\| \phi \|^{2}_{H^{1}_{\omega}(\mrI)} \ = \ 
  \sum_{i = 0}^{\infty} (1 + i^{2})^{1} \, \phi_{i}^{2} \, < \, \infty \, .
  \label{jjh1}
\ee

Using \eqref{propMnx},
\be
D \omega \phi(x) \ = \ - \omega^{(\alpha - \beta - 1 \, , \, \beta - 1)} \, 
  \sum_{i = 0}^{\infty} (1 + i) \, \phi_{i} \, 
  \frac{ \| | G_{i + 1}^{(\alpha - \beta - 1 \, , \, \beta - 1)} | \|}%
{ \| | G_{i}^{(\alpha - \beta \, , \,  \beta)} | \|} \, \what{G}_{i + 1}^{(\alpha - \beta - 1 \, , \, \beta - 1)} (x) 
\ := \  - \omega^{(\alpha - \beta - 1 \, , \, \beta - 1)} \, \eta(x) \, .
\label{jjh4}
\ee
Next, using \eqref{GnNormRelation},
\be
\| \eta \|^{2}_{L^{2}_{\omega^{(\alpha - \beta - 1 \, , \, \beta - 1)}}} \ \simeq \
\sum_{i = 0}^{\infty} ( (1 + i) \, \phi_{i} )^{2} \ \simeq \ \sum_{i = 0}^{\infty} (1 + i^{2}) \, \phi_{i}^{2}
\ \simeq  \ \| \phi \|^{2}_{H^{1}_{\omega}(\mrI)} \, .
\label{jjh5}
\ee


For $\psi(x) \, = \, \sum_{j = 0}^{\infty} \psi_{j} \, \what{G}_{j}^{(\beta , \alpha - \beta)}(x) \, \in 
 H^{\alpha - 1}_{\omega^{*}}(\mrI)$,
\be
\| \psi \|^{2}_{H^{\alpha - 1}_{\omega^{*}}(\mrI)} \ = \ 
  \sum_{j = 0}^{\infty} (1 + j^{2})^{\alpha - 1} \, \psi_{j}^{2} \, < \, \infty \, .
  \label{jjh3}
\ee

Using \eqref{propmcM},
\[
\mcM^{*} \omega^{*} \psi(x) \ = \ \sum_{j = 0}^{\infty} \mu_{j} \, \, \psi_{j} \, 
\frac{ \| | G_{j + 1}^{(\alpha - \beta - 1 \, , \, \beta - 1)} | \|}%
{ \| | G_{j}^{(\beta \, , \, \alpha - \beta)} | \|} \, \what{G}_{j + 1}^{(\alpha - \beta - 1 \, , \, \beta - 1)} (x) \, .
\]
Again using \eqref{GnNormRelation},
\be
\| \mcM^{*} \omega^{*} \psi \|^{2}_{L^{2}_{\omega^{(\alpha - \beta - 1 \, , \, \beta - 1)}}} \ \simeq \
\sum_{j = 0}^{\infty} ( \mu_{j} \, \psi_{j} )^{2} \ \simeq \ \sum_{j = 0}^{\infty} (1 + j^{2})^{\alpha - 1} \, \psi_{j}^{2}
\ \simeq  \ \| \psi \|^{2}_{H^{\alpha - 1}_{\omega}(\mrI)} \, .
\label{jjh2}
\ee

Combining the above pieces we have that
\begin{align}
|B_{0}(\phi , \psi)| &= | \left( k(x) \, D \omega \phi(x) \, , \, \mcM^{*} \omega^{*} \psi(x) \right) | \nonumber \\
&= | \left( - \, k(x) \, \eta(x) \, , \, \mcM^{*} \omega^{*} \psi(x) \right)_{\omega^{(\alpha - \beta - 1 \, , \, \beta - 1)}} |  \ 
\mbox{ (using \eqref{jjh4})} \nonumber \\
&\le \ \| k \|_{L^{\infty}(\mrI)}  \, 
      \| \eta \|_{L^{2}_{\omega^{(\alpha - \beta - 1 \, , \, \beta - 1)}}}  
      \, \| \mcM^{*} \omega^{*} \psi \|_{L^{2}_{\omega^{(\alpha - \beta - 1 \, , \, \beta - 1)}}} \nonumber  \\
&\lesssim \ \| k \|_{L^{\infty}(\mrI)}  \, \| \phi \|_{H^{1}_{\omega^{*}}(\mrI)} \,  \| \psi \|_{H^{\alpha - 1}_{\omega}(\mrI)} 
    \mbox{ (using \eqref{jjh5} and \eqref{jjh2})} .  \label{jjh6}
\end{align}    
 
\underline{Continuity of $B_{1}(\cdot , \cdot)$} \\
\label{ssec_B2}

For $\phi \in H^{1}_{\omega}(\mrI)$, an application of Theorem 5.1 in \cite{erv191} 
(see also the Proof of Lemma 3.1 in \cite{zhe211}) establishes  that $\omega \phi \, \in 
H^{2 - \alpha}_{(\beta - 1 \, , \, \alpha - \beta - 1)}(\mrI)$ with
$\| \omega \phi \|_{H^{2 - \alpha}_{(\beta - 1 \, , \, \alpha - \beta - 1)}(\mrI)} \, 
\lesssim \, \| \phi \|_{H^{1}_{\omega}(\mrI)}$.
Then, using Lemma~\ref{lmamapD}, $D \omega \phi \, \in 
H^{-(\alpha - 1)}_{\omega^{*}}(\mrI)$ with
$\| D \omega \phi \|_{H^{-(\alpha - 1)}_{\omega^{*}}(\mrI)} \, 
\lesssim \, \| \omega \phi \|_{H^{2 - \alpha}_{(\beta - 1 \, , \, \alpha - \beta - 1)}(\mrI)}  \, 
\lesssim \, \| \phi \|_{H^{1}_{\omega}(\mrI)}$.

Thus,
\begin{align}
|B_{1}(\phi , \psi)| &= | \langle b \, D u  \, , \, \psi \rangle_{\omega^{*}}|    \nonumber \\
&\le \ \| D \omega \phi  \|_{H^{-(\alpha - 1)}_{\omega^{*}}(\mrI)}  \, 
  \| b \,  \psi  \|_{H^{\alpha - 1}_{\omega^{*}}(\mrI)}   \nonumber \\
&\lesssim \  \| b \|_{W^{1, \infty}_{w}(I)} \, \| \phi \|_{H^{1}_{\omega}(\mrI)}
\|  \psi  \|_{H^{\alpha - 1}_{\omega^{*}}(\mrI)}  
  \, .  \label{jjh8} 
\end{align}

\underline{Continuity of $B_{2}(\cdot , \cdot)$} \\
To establish the continuity of 
$B_{2}(\cdot , \cdot)$, note that
\begin{align}
|B_{2}(\phi , \psi)| &= | \langle c u  \, , \, \psi \rangle_{\omega^{*}} |\
  = \bigg| \int_{I} \omega^{*}(x) \, c(x) \, \omega(x) \phi(x) \,  \psi(x) \, dx  \bigg|\nonumber \\
&\le \  \| c \|_{L^{\infty}(I)} \, \|  \omega^{1/2} \,  \omega^{* 1/2} \|_{L^{\infty}(\mrI)}
 \int_{I}  | \omega^{1/2}(x) \phi(x) | \,  | \omega^{* 1/2}(x) \, \psi(x) | \, dx \nonumber  \\
&\le \   \| c \|_{L^{\infty}(I)} \, \| \phi \|_{L^{2}_{\omega}(\mrI)} \, \| \psi \|_{L^{2}_{\omega^{*}}(\mrI)} \nonumber  \\
&\le \  \| c \|_{L^{\infty}(I)} \, \| \phi \|_{H^{1}_{\omega}(\mrI)} \, \| \psi \|_{H^{\alpha - 1}_{\omega^{*}}(\mrI)}   \label{jjh7}. 
\end{align}

Combining \eqref{jjh6}, \eqref{jjh8}, and \eqref{jjh7} we obtain the following lemma.
\begin{lemma} \label{lmaBcts}
For $k \in L^{\infty}(\mrI)$, $b \in W^{1 , \infty}_{w}(\mrI)$ and 
$c \in L^{\infty}(\mrI)$, there exists a constant $C_{1} > 0$ such that
$| B(\phi , \psi) | \, \le \, C_{1} \, \| \phi \|_{H^{1}_{\omega}(\mrI)} \, \| \psi \|_{H^{\alpha - 1}_{\omega^{*}}(\mrI)}$  for all
 $\phi \in H^{1}_{\omega}(\mrI)$ and   $ \psi \in H^{\alpha - 1}_{\omega^{*}}(\mrI)$.
 \end{lemma} 
\mbox{ } \hfill \qed


\subsection{Property \eqref{BB2} -- Positivity of $B(\phi , \cdot)$}
\label{ssec_Pos}

To establish the positivity of $B(\phi , \cdot)$, i.e. property \eqref{BB2}, we need to show that there exists $C_{2} > 0$
such that for
$\phi \in H^{1}_{\omega}(\mrI)$
\[
  \sup_{\psi \in H^{1}_{\omega^{*}}(\mrI)}  B(\phi , \psi)  \ \ge \ C_{2} \,  \| \phi \|_{H^{\alpha - 1}_{\omega}(\mrI)} \,
   \| \psi \|_{H^{\alpha - 1}_{\omega^{*}}(\mrI)} \,  .
\]

In comparison with the case $k = 1$ discussed in \cite{zhe211}, given $\phi \in H^{1}_{\omega}(\mrI)$ the appropriate
choice of $\psi \in H^{\alpha - 1}_{\omega^{*}}(\mrI)$ for $k(x) \neq \, constant$ is more restrictive. Specifically, in this
case we need $\psi$ such that 
$$D (\omega \phi(x) ) \mcM^{*} (\omega^{*} \psi(x) )> 0$$ for $x \, a.e.$ in $\mrI$. 
This then 
allows $k(x)$ to be taken out of the inner product and the orthogonality property of the Jacobi polynomials to be used.

To achieve this we choose $\phi$ and $\psi$ such that $D (\omega \phi(x))  = \mcM^{*} (\omega^{*} \psi(x))$. Let 
$$\phi(x) = \sum_{i = 0}^{\infty} \phi_{i}  \what{G}_{i}^{(\alpha - \beta \, , \, \beta)} (x)\in H^1_\omega(\mrI),$$
 and define $\Phi_i:=-\phi_i/\mu_i$ such that
\begin{align}
\phi(x) &= \ - \, \sum_{i = 0}^{\infty} \mu_{i}  \, \Phi_{i} \, \what{G}_{i}^{(\alpha - \beta \, , \, \beta)} (x) , 
 \label{defphi1}  \\
 \mbox{and } \ 
\psi(x) &=  \  \sum_{j = 0}^{\infty} (j + 1) \, \Phi_{j} \, \what{G}_{j}^{(\beta \, , \, \alpha - \beta)} (x) \, .
\label{defpsi1}
\end{align}

As $\phi \in H^{1}_{\omega}(\mrI)$, using \eqref{propmcM2},
\[
\| \phi \|_{H^{1}_{\omega}(\mrI)}^{2} \ = \ 
\sum_{i= 0}^{\infty} (1 + i^{2})^{1} \left( \mu_{i} \, \Phi_{i} \right)^{2} \ \simeq \ 
\sum_{i= 0}^{\infty} (1 + i^{2})^{1} \, (1 + i)^{2 ( \alpha - 1)} \, \Phi_{i}^{2} \ \simeq \ 
\sum_{i= 0}^{\infty} (1 + i^{2})^{\alpha} \, \Phi_{i}^{2} \, < \, \infty \, .
\]

Then,
\[
\| \psi \|_{H^{\alpha - 1}_{\omega^{*}}(\mrI)}^{2} \ = \ 
\sum_{j= 0}^{\infty} (1 + j^{2})^{\alpha - 1} \left( (j + 1) \, \Phi_{j} \right)^{2} \ \simeq \ 
\sum_{j= 0}^{\infty} (1 + j^{2})^{\alpha} \, \Phi_{j}^{2} \, < \, \infty \, .
\]
i.e., $\psi \in H^{\alpha - 1}_{\omega^{*}}(\mrI)$ and 
$ \| \psi \|_{H^{\alpha - 1}_{\omega^{*}}(\mrI)} \simeq \| \phi \|_{H^{1}_{\omega}(\mrI)}^{2}$.

Using \eqref{propmcM} and \eqref{propMnx}
\begin{align*}
B_{0}(\phi , \psi) &= \ \left( k(x) ( D \omega \phi(x) ) ,  \mcM^{*} (\omega^{*} \psi(x)) \right)  \nonumber \\
&= \ \left( k(x) \, \omega^{(\alpha - \beta - 1 \, , \, \beta - 1)} \,
\sum_{i = 0}^{\infty} (i + 1) \, \mu_{i} \, \Phi_{i} \, \frac{ \| | G_{i + 1}^{(\alpha - \beta - 1 \, , \, \beta - 1)} | \|}%
{ \| | G_{i}^{(\alpha - \beta \, , \, \beta)} | \|} \, \what{G}_{i + 1}^{(\alpha -\beta - 1 \, , \,  \beta - 1)} (x) \ , \ \right. \\
& \quad \quad \quad \quad  \quad \quad \quad \quad   \quad \quad \quad   \left. 
\sum_{j = 0}^{\infty} (j + 1) \, \mu_{j} \, \Phi_{j} \, \frac{ \| | G_{j + 1}^{(\alpha - \beta - 1 \, , \, \beta - 1)} | \|}%
{ \| | G_{j}^{(\beta \, , \, \alpha - \beta)} | \|} \, \what{G}_{j + 1}^{(\alpha - \beta - 1 \, , \, \beta - 1)} (x) \right) \, .
\end{align*}

As $k(x)$ and $\omega^{(\alpha - \beta - 1 \, , \, \beta - 1)}(x) \ge 0$ for all $x \in \mrI$, 
$\| | G_{i}^{(\alpha - \beta \, , \, \beta)} | \| \, = \, \| | G_{i}^{(\beta \, , \, \alpha - \beta)} | \|$,
and using \eqref{GnNormRelation},
\begin{align}
B_{0}(\phi , \psi) &\gtrsim \ k_{0} \ 
\left( \omega^{(\alpha - \beta - 1 \, , \, \beta - 1)} \, 
\sum_{i = 0}^{\infty} (i + 1) \, \mu_{i} \, \Phi_{i}  \, \what{G}_{i + 1}^{(\alpha - \beta - 1 \, , \, \beta - 1)} (x) \ ,  \
  \sum_{j = 0}^{\infty} (j + 1) \, \mu_{j} \, \Phi_{j} \,  \what{G}_{j + 1}^{(\alpha - \beta - 1 \, , \, \beta - 1)} (x) \right)  
  \nonumber  \\
&= \ k_{0} \, \sum_{i = 0}^{\infty} (i + 1)^{2} \, \mu_{i}^{2} \, \Phi_{i}^{2} 
 \ \simeq \ k_{0} \,  \sum_{i = 0}^{\infty} (1 + i^{2}) \, \phi_{i}^{2}  
\ \simeq  \ k_{0} \, \| \phi \|_{H^{1}_{\omega}(\mrI)}^{2}   \nonumber \\
&\simeq  \ k_{0} \, \| \phi \|_{H^{1}_{\omega}(\mrI)} \, \| \psi \|_{H^{\alpha - 1}_{\omega^{*}}(\mrI)} \, . 
 \label{coeB0_1}
\end{align} 

Combining \eqref{coeB0_1}, \eqref{jjh8}, and \eqref{jjh7} we obtain the following lemma.
\begin{lemma} \label{lmaBpos}
For $0 < k_{0} \le k(x)  \in L^{\infty}(\mrI)$, and $b \in W^{1 , \infty}_{w}(\mrI)$, 
$c \in L^{\infty}(\mrI)$ sufficiently small (with respect to $k_{0}$),
there exists a constant $C_{2} > 0$ such that
\be
  \sup_{\psi \in H^{\alpha - 1}_{\omega^{*}}(\mrI)}  B(\phi , \psi)  \ \ge \ C_{2} \,  \| \phi \|_{H^{1}_{\omega}(\mrI)} \,
   \| \psi \|_{H^{\alpha - 1}_{\omega^{*}}(\mrI)}  \ \mbox{ for all } \phi  \in H^{1}_{\omega}(\mrI)  \, .
\label{Bposv1}
\ee
 \end{lemma} 
\textbf{Proof}:
Let $\phi  \in H^{1}_{\omega}(\mrI)$ be given by \eqref{defphi1} and choose 
$\psi \in H^{\alpha - 1}_{\omega^{*}}(\mrI)$ to be given by \eqref{defpsi1}. 
Then, from \eqref{coeB0_1}, \eqref{jjh8}, and \eqref{jjh7}
there exists $c_{0}, \, c_{1} > 0$ such that
\begin{align*}
B(\phi , \psi) &\ge \ B_{0}(\phi , \psi) \ - \ | B_{1}(\phi , \psi)  \ - \ B_{2}(\phi , \psi) |  \\
&\ge \ c_{0} \, k_{0} \, \| \phi \|_{H^{1}_{\omega}(\mrI)} \, \| \psi \|_{H^{\alpha - 1}_{\omega^{*}}(\mrI)}
 \ - \ c_{1}  \| b \|_{W^{1, \infty}_{w}(I)} \, \| \phi \|_{H^{1}_{\omega}(\mrI)}
\| \psi  \|_{H^{\alpha - 1}_{\omega^{*}}(\mrI)}     \\
& \quad \quad \quad \quad  \quad \quad \quad 
\ - \ \| c \|_{L^{\infty}(I)} \, \| \phi \|_{H^{1}_{\omega}(\mrI)} \, \| \psi \|_{H^{\alpha - 1}_{\omega^{*}}(\mrI)}  \\
&\ge \ C_{2} \, \| \phi \|_{H^{1}_{\omega}(\mrI)} \, \| \psi \|_{H^{\alpha - 1}_{\omega^{*}}(\mrI)} \, , \
\end{align*}
for $b \in W^{1 , \infty}_{w}(\mrI)$ and $c \in L^{\infty}(\mrI)$  sufficiently small.  
The stated result \eqref{Bposv1} then follows. \\
\mbox{ } \hfill \qed

\subsection{Property \eqref{BB3} -- Positivity of $B(\cdot , \psi)$}
\label{ssec_Pos2}

The argument to establish \eqref{BB3} is almost verbatim of that used to show \eqref{BB2}.
\begin{lemma} \label{lmaBpos2}
For $0 < k_{0} \le k(x)  \in L^{\infty}(\mrI)$, and $b \in W^{1 , \infty}_{w}(\mrI)$, 
$c \in L^{\infty}(\mrI)$ sufficiently small (with respect to $k_{0}$),
there exists a constant $C_{2} > 0$ such that
\be
  \sup_{\phi \in H^{1}_{\omega}(\mrI)}  B(\phi , \psi)  \ \ge \ C_{2} \,  \| \phi \|_{H^{1}_{\omega}(\mrI)} \,
   \| \psi \|_{H^{\alpha - 1}_{\omega^{*}}(\mrI)}  \ \mbox{ for all } \psi \in H^{\alpha - 1}_{\omega^{*}}(\mrI)  \, .
\label{Bposv2}
\ee
 \end{lemma} 
\textbf{Proof}:
Let $\psi  \in H^{1}_{\omega^{*}}(\mrI)$ be given by \eqref{defpsi1} and choose 
$\phi \in H^{\alpha - 1}_{\omega}(\mrI)$ to be given by \eqref{defphi1}. 

Then,  as in the proof of Lemma \ref{lmaBpos}, we obtain
$B(\phi , \psi) \ \ge \ C_{2} \, \| \phi \|_{H^{1}_{\omega}(\mrI)} \, \| \psi \|_{H^{\alpha - 1}_{\omega^{*}}(\mrI)} $,
for $b \in W^{1 , \infty}_{w}(\mrI)$ and $c \in L^{\infty}(\mrI)$  sufficiently small.   \\
\mbox{ } \hfill \qed

Combining the above we obtain the following existence and uniqueness result for the solution of \eqref{wform1}.
\begin{theorem} \label{thmexun1}
For $f \in H^{-(\alpha - 1)}_{\omega^{*}}(\mrI)$, $0 < k_{0} \le k(x)  \in L^{\infty}(\mrI)$, and $b \in W^{1 , \infty}_{w}(\mrI)$, 
$c \in L^{\infty}(\mrI)$ sufficiently small (with respect to $k_{0}$),
there exists a unique  $u(x) \ = \ \omega \, \phi(x)$ satisfying \eqref{wform1} with $\phi  \in H^{1}_{\omega}(\mrI)$,
and  $\| \phi \|_{H^{1}_{\omega}(\mrI)} \, \le \, \frac{1}{C_{2}} \| f \|_{H^{-(\alpha - 1)}_{\omega^{*}}(\mrI)}$. 
 \end{theorem} 

\textbf{Proof}: Note that $F$ defined by \eqref{defF} satisfies
\begin{align*}
\| F \| &= \ \sup_{0 \ne \psi \in H^{\alpha - 1}_{\omega^{*}}(\mrI)} \frac{ | F(\psi) |}{ \| \psi \|_{H^{\alpha - 1}_{\omega^{*}}(\mrI)}}
\ = \ \sup_{0 \ne \psi \in H^{\alpha - 1}_{\omega^{*}}(\mrI)}
 \frac{ | \langle f , \psi \rangle_{\omega^{*}} |}{ \| \psi \|_{H^{\alpha - 1}_{\omega^{*}}(\mrI)}}    \\
&\le \ \sup_{0 \ne \psi \in H^{\alpha - 1}_{\omega^{*}}(\mrI)}
 \frac{ \| f \|_{H^{-(\alpha - 1)}_{\omega^{*}}(\mrI)} \| \psi \|_{H^{\alpha - 1}_{\omega^{*}}(\mrI)} }%
 { \| \psi \|_{H^{\alpha - 1}_{\omega^{*}}(\mrI)}}
 \ = \ \| f \|_{H^{-(\alpha - 1)}_{\omega^{*}}(\mrI)} \, .
 \end{align*}
Hence, $F$ defines a bounded linear functional.
The existence, uniqueness and bound for $\phi$ then follows from 
combining Lemmas \ref{lmaBcts}, \ref{lmaBpos} and \ref{lmaBpos2} with Theorem \ref{BNBthm}. 
\mbox{ } \hfill \qed

\subsection{The case of $\grave{\mcL}^{\alpha}_{r} \cdot$}
\label{ssec_Odo1}
For model \eqref{DefProb2} with $\wtilde{\mcL}_{r}^{\alpha}=\grave{\mcL}_{r}^{\alpha}$ we have the following problem: 
Given $f \in H^{-1}_{\omega^{*}}(\mrI)$,  $0 < k_{0} \le k(x)  \in L^{\infty}(\mrI)$,
$b \in W^{1 , \infty}_{w}(\mrI)$ and 
$c \in L^{\infty}(\mrI)$, 
determine $\phi \in H^{\alpha - 1}_{\omega}(\mrI)$
such that $u(x) \ = \ \omega(x) \, \phi(x)$ satisfies
\be
\langle \grave{\mcL}_{r}^{\alpha} u \ + \ b \, D u \ +  \ c \, u \, , \, \psi \rangle_{\omega^{*}} 
\ = \ \langle f \, , \, \psi \rangle_{\omega^{*}} \, , \ \ \forall \, \psi \in H^{1}_{\omega^{*}}(\mrI) \, .
\label{wform1*}
\ee

The analysis for the existence and uniqueness of solution to \eqref{wform1*} follows similarly to that presented
for \eqref{wform1}. In place of \eqref{propmcM} the following result is used.
From \cite{jia181},  with $\mcM u(x)\, := \, \big( r D^{-(2 - \alpha)} \ + \ (1 - r) D^{-(2 - \alpha)*}  \big) D u(x)$ 
\begin{align}
   \mcM  \omega(x) \, G_{k}^{(\alpha - \beta \, , \, \beta)}(x)
   &= \ \mu_{k} \, G_{k+1}^{(\beta - 1 \, , \, \alpha - \beta - 1)}(x) \, ,  \label{propmCM}  \\
   \mbox{ where } \ \ 
   \mu_{k} &= \ c_{*}^{*}
    \frac{\Gamma(k + \alpha)}{\Gamma(k + 1)}  \, ,
    \ k = 0, 1, 2, \ldots, \, .
\label{propmCM1}  
\end{align}

In summary we have the following result.
\begin{theorem} \label{thmexun2}
For $f \in H^{-1}_{\omega^{*}}(\mrI)$, $0 < k_{0} \le k(x)  \in L^{\infty}(\mrI)$, and $b \in W^{1 , \infty}_{w}(\mrI)$, 
$c \in L^{\infty}(\mrI)$ sufficiently small (with respect to $k_{0}$),
there exists a unique  $u(x) \ = \ \omega \, \phi(x)$ satisfying \eqref{wform1*} with $\phi  \in H^{\alpha - 1}_{\omega}(\mrI)$,
and  $\| \phi \|_{H^{\alpha - 1}_{\omega}(\mrI)} \, \le \, \frac{1}{C_{2}} \| f \|_{H^{-1}_{\omega^{*}}(\mrI)}$. 
\end{theorem}

 \setcounter{equation}{0}
\setcounter{figure}{0}
\setcounter{table}{0}
\setcounter{theorem}{0}
\setcounter{lemma}{0}
\setcounter{corollary}{0}
\setcounter{definition}{0}

\section{Regularity of the solutions}
\label{sec_reg}
In this section we investigate the regularity of the solution to \eqref{DefProb2}-\eqref{DefBC2}.
The first step is to characterize the mapping properties of 
$( \wtilde{\mcL}_{r}^{\alpha} )^{-1}$.
With this in hand the second step employs a boot strapping argument to determine the regularity of the solution. We again give the detailed analysis for $\acute{\mcL}_{r}^{\alpha} t$ and then give the corresponding results
for $\grave{\mcL}_{r}^{\alpha} \cdot$.

\subsection{Mapping properties of $\acute{\mcL}_{r}^{\alpha} \cdot$}
\label{ssec_map1}
The following two lemmas are helpful in determining the mapping  properties of $\acute{\mcL}_{r}^{\alpha} \cdot$.
We firstly introduce the subspace of $H^{t}_{(a , b)}(\mrI)$, $_{0}H^{t}_{(a , b)}(\mrI)$, defined by
\[
  _{0}H^{t}_{(a , b)}(\mrI) \ := \ \{ h \in H^{t}_{(a , b)}(\mrI) \, : \, \langle h \, , \, 1 \rangle_{\omega^{(a , b)}} = 0 \} \, .
\]
In words, $_{0}H^{t}_{(a , b)}(\mrI) $ represents those functions in $H^{t}_{(a , b)}(\mrI)$ that when expanded
in the basis $ \{ \what{G}_{i}^{(a , b)}(x) \}_{i = 0}^{\infty}$ have coefficient of $\what{G}_{0}^{(a , b)}(x)$ equal to 0.

\begin{lemma} \label{lmaregN}
For $1 < \alpha < 2$, $0 \le r \le 1$, let $\beta$ be determined be \textbf{Condition A} and let
$\mcN \, : \, \omega^{(\alpha - \beta - 1 \, , \, \beta - 1)}(x) \otimes  \mbox{}_{0}H^{s}_{(\alpha - \beta - 1 \, , \, \beta - 1)}(\mrI)
\ \rightarrow \ H^{s - (\alpha - 1)}_{(\beta  , \, \alpha - \beta)}(\mrI)$ be defined by
\be
  \mcN \, \omega^{(\alpha - \beta - 1 \, , \, \beta - 1)}(x) \, \psi(x) \ := \ 
  - D \big( r D^{-(2 - \alpha)} \ + \ (1 - r) D^{-(2 - \alpha)*}  \big) \omega^{(\alpha - \beta - 1 \, , \, \beta - 1)}(x) \, \psi(x) \, .
\label{defmcN}
\ee
Then the mapping $\mcN$  is bijective, continuous and has a continuous inverse.
\end{lemma}
\textbf{Proof}: We begin by showing that $\mcN$ is well defined.

From \cite[Lemma 2.2]{jia181}, we have for $k = 0, 1, 2, \ldots$
\begin{align}
\big( r D^{-(2 - \alpha)} \ + \ (1 - r) D^{-(2 - \alpha)*}  \big) \omega^{(\alpha - \beta - 1 \, , \, \beta - 1)}(x) 
G_{k}^{(\alpha - \beta - 1 \, , \, \beta - 1)}(x) &= \ \sigma_{k} \, G_{k}^{(\beta - 1 \, , \, \alpha - \beta - 1)}(x) 
\, ,     \label{NappG}  \\
\mbox{where } \ \ \sigma_{k} \ := \ - c_{*}^{*} \frac{\Gamma(k + \alpha -1)}{\Gamma(k + 1)} \, . & 
\label{defSig}
\end{align}

Thus, using \eqref{eqC4},
\[
 D \big( r D^{-(2 - \alpha)} \ + \ (1 - r) D^{-(2 - \alpha)*}  \big) \omega^{(\alpha - \beta - 1 \, , \, \beta - 1)}(x) 
G_{k}^{(\alpha - \beta - 1 \, , \, \beta - 1)}(x) \ = \ - \mu_{k} \, G_{k - 1}^{(\beta  , \, \alpha - \beta)}(x) 
\, , \ k = 1, 2, \ldots ,   
\]
where $\mu_{k}$ is given by \eqref{propmcM1}.

Now, for $\psi(x) \ = \ \sum_{i = 1}^{\infty} \psi_{i} \, \what{G}_{i}^{(\alpha - \beta - 1 \, , \, \beta - 1)}(x)
\in \mbox{}_{0}H^{s}_{(\alpha - \beta - 1 \, , \, \beta - 1)}(\mrI)$,
\[
\| \psi \|_{H^{s}_{(\alpha - \beta - 1 \, , \, \beta - 1)}(\mrI)}^{2} \ = \ 
\sum_{i = 1}^{\infty} (1 + i^{2})^{s} \, \psi_{i}^{2} \ < \ \infty \, .
\]

Note that $\mcN \omega^{(\alpha - \beta - 1 \, , \, \beta - 1)}(x) \, \psi(x) \ = \ 
  \sum_{i = 1}^{\infty} \mu_{i} \, \psi_{i} \,
 \frac{ |\| G_{i - 1}^{(\beta  , \, \alpha - \beta)} \| |}{  |\| G_{i}^{(\alpha - \beta - 1  , \, \beta - 1)} \| |}
 \what{G}_{i - 1}^{(\beta  , \alpha - \beta )}(x) $, and using \eqref{GnNormRelation}  and \eqref{propmcM2},
\begin{align*}
\| \mcN \omega^{(\alpha - \beta - 1 \, , \, \beta - 1)}(x) \, \psi(x) \|_{H^{s - (\alpha - 1)}_{(\alpha - \beta - 1 \, , \, \beta - 1)}(\mrI)}^{2} 
&\simeq \ \sum_{i = 1}^{\infty} (1 + (i - 1)^{2})^{s - (\alpha - 1)} \,  \left( \mu_{i} \, \psi_{i} \right)^{2}  \\
&\simeq \ \sum_{i = 1}^{\infty} (1 + i^{2})^{s - (\alpha - 1)} \,  (1 + i)^{2 (\alpha - 1) } \, \psi_{i}^{2}   \\
&\simeq \ \sum_{i = 1}^{\infty} (1 + i^{2})^{s} \, \psi_{i}^{2}  \ 
\simeq \ \| \psi \|_{H^{s}_{(\alpha - \beta - 1 \, , \, \beta - 1)}(\mrI)}^{2} \ < \ \infty \, .
\end{align*}

Hence it follows that $\mcN$  is a continuous, bijective mapping from 
$\omega^{(\alpha - \beta - 1 \, , \, \beta - 1)}(x) \otimes  \mbox{}_{0}H^{s}_{(\alpha - \beta - 1 \, , \, \beta - 1)}(\mrI)$
onto $H^{s - (\alpha - 1)}_{(\beta  , \, \alpha - \beta }(\mrI)$ with a continuous inverse. \\
\mbox{ } \hfill \qed

The second lemma involves the differentiation operator, $D$, and the identity given by \eqref{eqC4}.
\begin{lemma} \label{lmaregD}
Let $a, b > 0$. Then,
\be
  D \, : \, \omega^{(a , b)}(x) \otimes H_{(a , b)}^{t}(\mrI) \ \rightarrow \ 
               \omega^{(a - 1 \, , \, b - 1)}(x) \otimes  \mbox{}_{0}H_{(a - 1 \, , \, b - 1)}^{t - 1}(\mrI)  
  \label{Dmapr1}
\ee
is bijective, continuous and has a continuous inverse.
\end{lemma}
\textbf{Proof}: We begin by showing that $D$, with domain and range space specified in \eqref{Dmapr1} is well defined. 

Let $\phi(x) \ = \ \sum_{i = 0}^{\infty} \phi_{i} \, \what{G}_{i}^{(a , b)}(x) \ \in \ H_{(a , b)}^{t}(\mrI)$. Then,
\[
 \| \phi \|^{2}_{H_{(a , b)}^{t}(\mrI)} \ = \ \sum_{i = 0}^{\infty} (1 + i^{2})^{t} \, \phi_{i}^{2} \ < \ \infty \, .
\]

Using  \eqref{eqC4p5},
\begin{align*}
D \omega^{(a , b)}(x) \, \phi(x) &= \ 
   \omega^{(a - 1 \, , \, b - 1)}(x) \, \sum_{i = 0}^{\infty} - (i + 1) \, \phi_{i} \, 
  \frac{ |\| G_{i + 1}^{(a - 1 \, , \, b - 1)} \| |}{  |\| G_{i}^{(a  ,  b)} \| |}
 \what{G}_{i + 1}^{(a - 1 \, , \, b - 1)}(x)   \\
&:= \ \omega^{(a - 1 \, , \, b - 1)}(x) \, \psi(x) \, .
\end{align*}

Then, using \eqref{GnNormRelation},
\begin{align*}
\| \psi \|^{2}_{H_{(a - 1 \, , \, b - 1)}^{t - 1}(\mrI)} 
&\simeq \ \sum_{i = 0}^{\infty} (1 + i^{2})^{t - 1} \, ( (i + 1) \, \phi_{i} )^{2}   \\
&\simeq \ \sum_{i = 0}^{\infty} (1 + i^{2})^{t} \, \phi_{i}^{2} \ \simeq \ \| \phi \|^{2}_{H_{(a , b)}^{t}(\mrI)}  \, .
\end{align*}

Hence it follows that $D$  is a continuous, bijective mapping from 
$\omega^{(a , b)}(x) \otimes H_{(a , b)}^{t}(\mrI)$
onto $\omega^{(a - 1 \, , \, b - 1)}(x) \otimes  \mbox{}_{0}H_{(a - 1 \, , \, b - 1)}^{t - 1}(\mrI)$ with a continuous inverse. \\
\mbox{ } \hfill \qed

Recall, by the imbedding of the Sobolev spaces, $H^{s}_{(a , b)}(\mrI) \subset H^{s'}_{(a , b)}(\mrI)$ for $s' < s$.

For $f \in H^{s}_{(\beta , \alpha - \beta)}(\mrI)$, $s \ge -1$, from Theorem \ref{thmexun1} we have that for
$0 < k_{0} \le k(x)  \in L^{\infty}(\mrI)$ there exists a unique solution $u(x) \ = \ \omega \, \phi(x)$ satisfying
\[
   \acute{\mcL}_{r}^{\alpha} \omega \, \phi(x) \ = \ \mcN \, k(x) \, D \omega \, \phi(x) \ = \ f(x) \, .
\]

As $f \in H^{s}_{(\beta , \alpha - \beta)}(\mrI)$, by the mapping properties of $\mcN$, there exists a unique 
$\psi \in  \mbox{}_{0}H^{s + (\alpha - 1)}_{(\alpha - \beta - 1 \, , \, \beta - 1)}(\mrI)$ such that
\[
  k(x) \, D \omega \, \phi(x) \ = \ \omega^{(\alpha - \beta - 1 \, , \, \beta - 1)}(x) \, \psi(x) \, .
\]

Assuming that $\phi(x) \in H^{t}_{(\alpha - \beta , \beta)}(\mrI)$ for some $t \, \ge \, \alpha - 1$, from the mapping
properties of $D$ given by Lemma \ref{lmaregD}, there exists 
$\wtilde{\phi}(x) \in \mbox{}_{0}H^{t - 1}_{(\alpha - \beta - 1 \, , \, \beta - 1)}(\mrI)$
such that
\[
    k(x) \, \wtilde{\phi}(x) \ = \ \psi(x) \,  \in  \, H^{s' + (\alpha - 1)}_{(\alpha - \beta - 1 \, , \, \beta - 1)}(\mrI)  \ \mbox{ for } \
   s'  \le  s  \, .
\]

The regularity of $k(x) \, \wtilde{\phi}(x)$ and $\psi(x) \in  \, H^{s' + (\alpha - 1)}_{(\alpha - \beta - 1 \, , \, \beta - 1)}(\mrI) $
must match. Specifically, if $k(x)$ is not sufficiently smooth $\psi(x)$ must be interpreted to line in some
$H^{s' + (\alpha - 1)}_{(\alpha - \beta - 1 \, , \, \beta - 1)}(\mrI)$ space for $s' < s$ in order to match the regularity
of $k(x) \, \wtilde{\phi}(x)$.

In view of Lemma \ref{lmaprodsp} and the above discussion we have the following theorem.
\begin{theorem}  \label{thexr1}
For $f \in H^{s}_{\omega^{*}}(\mrI)$, $-(\alpha - 1) \le s$, $0 < k_{0} \le k(x)  
\in W_{w}^{\ceil{s + \alpha - 1} , \infty}(\mrI)$, 
there exists a unique  $u(x) \ = \ \omega \, \phi(x)$ satisfying 
$\acute{\mcL}_{r}^{\alpha} \omega \, \phi(x) \ =  \ f(x)$  with $\phi  \in H^{s + \alpha}_{\omega}(\mrI)$.
 \end{theorem} 

As shown in \cite{erv191, hao201}, the presence of an advection or reaction term limits the regularity of the solution to 
\eqref{DefProb2}, \eqref{DefBC2}. The regularity of the solution of \eqref{DefProb2}, \eqref{DefBC2} for $k(x) = 1$
was given by Theorems 5.2 and 5.3 in \cite{erv191} (see also \cite{zhe211}). Combining that result with 
Theorems \ref{thmexun1} and \ref{thexr1} we  have the following corollary.

For $\epsilon > 0$ arbitrary, introduce $\wtilde{s}$ defined by
\be
\wtilde{s} \, := \, \left\{ \begin{array}{rl}
\min\{s , \, \alpha + (\alpha - \beta) + 1 - \epsilon, \, \alpha +  \beta + 1 - \epsilon \} , \ & \mbox{ if } b = 0, \ 
 \\
\min\{s , \, \alpha + (\alpha - \beta) - 1 - \epsilon, \, \alpha +  \beta - 1 - \epsilon \} , \ & \mbox{ if } b \ne 0  .
\end{array} \right.
\label{defsar}
\ee

\begin{corollary}  \label{corregk1}
For $f \in H^{s}_{\omega^{*}}(\mrI)$, $-(\alpha - 1) \le s$, $\beta$ be determined by \textbf{Condition A},
$0 < k_{0} \le k(x)  
\in W_{w}^{\ceil{s + \alpha - 1} , \infty}(\mrI)$, 
 $b(x) \in W^{\max\{1 , \, \ceil{ \, \wtilde{s} \,} \}, \infty}_{w}(\mrI)$
and $c(x) \in 
W^{\ceil{ \, \wtilde{s} \,} , \infty}_{w}(\mrI)$
sufficiently small (with respect to $k_{0}$),
there exists a unique  $u(x) \ = \ \omega \, \phi(x)$ satisfying \eqref{wform1} with 
$\phi  \in 
H^{\alpha \, + \,\wtilde{s}}_{(\alpha - \beta \, , \, \beta)}(\mrI)$. In addition, there exists $C > 0$ such that
\be
   \| \phi \|_{H^{\alpha \, + \, \wtilde{s}}_{(\alpha - \beta \, , \, \beta)}(\mrI)} \ \le \ C \,  \| f \|_{H^{\wtilde{s}}_{\omega^{*}}(\mrI)} \, .
\label{kjg1}
\ee
\end{corollary}

\textbf{Remark}: 
 The norm
estimate \eqref{kjg1} follows from that at each of the (finite number of) steps in the 
boot strapping argument used to establish the regularity of $\phi$ 
the terms on the right hand side are bounded by a constant times $ \| f \|_{H^{\wtilde{s}}_{\omega^{*}}(\mrI)}$. \\

\subsection{Mapping properties of $\grave{\mcL}_{r}^{\alpha} \cdot$}
\label{ssec_map2}

The mapping properties for $\grave{\mcL}_{r}^{\alpha} \cdot$ are obtained in a similar manner to 
$\acute{\mcL}_{r}^{\alpha} \cdot$. In place of Lemmas \ref{lmaregN} and \ref{lmaregD} we have the 
following two lemmas.

\begin{lemma} \label{lmaregD2}
Let $a, b > -1$. Then,
\be
  D \, : \, \mbox{}_{0}H_{(a  , \, b)}^{s}(\mrI)  \ \rightarrow \ 
              H_{(a + 1 \, , \, b + 1)}^{s - 1}(\mrI)  
  \label{Dmapr2}
\ee
is bijective, continuous and has a continuous inverse.
\end{lemma}
\textbf{Proof}: The proof uses property \eqref{eqC4}, and it similar to that for Lemma \ref{lmaregD}. \\
\mbox{ } \hfill \qed

\begin{lemma} \label{lmaregM}
For $1 < \alpha < 2$, $0 \le r \le 1$, let $\beta$ be determined be \textbf{Condition A} and let
$\mcM u(x)\, := \, \big( r D^{-(2 - \alpha)} \ + \ (1 - r) D^{-(2 - \alpha)*}  \big) D u(x)$. Then,
$\mcM \, : \, \omega^{(\alpha - \beta \, , \, \beta)}(x) \otimes  H^{s}_{(\alpha - \beta \, , \, \beta)}(\mrI)
\ \rightarrow \ \mbox{}_{0}H^{s - (\alpha - 1)}_{(\beta  - 1, \, \alpha - \beta - 1)}(\mrI)$ 
  is bijective, continuous and has a continuous inverse.
\end{lemma}
\textbf{Proof}: The proof uses \eqref{propmCM}, \eqref{propmCM1} and is similar to that for Lemma \ref{lmaregN}. \\
\mbox{ } \hfill \qed

Noting that $\grave{\mcL}_{r}^{\alpha} \omega \phi(x) \ = \ - D \, k(x) \, \mcM \omega \phi(x)$, a similar argument
as was used above for $\acute{\mcL}_{r}^{\alpha} \cdot$ results in the following theorem and corollary.

\begin{theorem}  \label{thexr2}
For $f \in H^{s}_{\omega^{*}}(\mrI)$, $-1 \le s$, $0 < k_{0} \le k(x)  
\in W_{w}^{\ceil{s + 1} , \infty}(\mrI)$, 
there exists a unique  $u(x) \ = \ \omega \, \phi(x)$ satisfying 
$\grave{\mcL}_{r}^{\alpha} \omega \, \phi(x) \ =  \ f(x)$  with $\phi  \in H^{s + \alpha}_{\omega}(\mrI)$.
 \end{theorem} 
 
 \begin{corollary}  \label{corregk1*}
For $f \in H^{s}_{\omega^{*}}(\mrI)$, $-1 \le s$, $\beta$ be determined by \textbf{Condition A},
$0 < k_{0} \le k(x)  
\in W_{w}^{\ceil{s + 1} , \infty}(\mrI)$, 
$b(x) \in W^{\max\{1 , \, \ceil{ \, \wtilde{s} \,} \}, \infty}_{w}(\mrI)$
and $c(x) \in 
W^{\ceil{ \, \wtilde{s} \,} , \infty}_{w}(\mrI)$
sufficiently small (with respect to $k_{0}$),
there exists a unique  $u(x) \ = \ \omega \, \phi(x)$ satisfying \eqref{wform1} with 
$\phi  \in 
H^{\alpha \, + \, \wtilde{s}}_{(\alpha - \beta \, , \, \beta)}(\mrI)$. In addition, there exists $C > 0$ such that
\be
   \| \phi \|_{H^{\alpha \, + \, \wtilde{s}}_{(\alpha - \beta \, , \, \beta)}(\mrI)} \ \le \ C \,  \| f \|_{H^{\wtilde{s}}_{\omega^{*}}(\mrI)} \, .
\label{kjg1*}
\ee
\end{corollary}

\setcounter{equation}{0}
\setcounter{figure}{0}
\setcounter{table}{0}
\setcounter{theorem}{0}
\setcounter{lemma}{0}
\setcounter{corollary}{0}
\setcounter{definition}{0}

\section{Numerical approximation and analysis}
\label{sec_appx}
In this section, we prove error estimates in both weighted $L^2$ and energy norms for a numerical approximation of (\ref{DefProb2})--(\ref{DefBC2}). Building on the analysis of the bi-linear form, 
$B(\cdot , \cdot)$, in Section \ref{sec_form},  we follow the ideas in \cite[Section 4]{zhe211} to perform numerical analysis for the different trial and test spaces used here. We emphasize that the novel analysis on the inf-sup condition for $B(\cdot , \cdot)$ in Section \ref{sec_form}, which resolves the loss of coercivity caused by the variable coefficient $k(x)$, is critical in the error estimates. This distinguishes the current work, by extending the analysis of  Petrov-Galerkin spectral approximation method of constant coefficient problems to the variable coefficient case.

As we have done above, we give the analysis for the diffusion operator $\acute{\mcL}_{r}^{\alpha}\cdot$ and then
the corresponding result for $\grave{\mcL}_{r}^{\alpha}\cdot$.

\subsection{Approximation of \eqref{wform1}}
\label{sec_appx1}
Let
$X_{N} \, := \, \mbox{span}  \{ \what{G}_{j}^{(\alpha - \beta \, , \, \beta)} \}_{j = 0}^{N} 
\subset H^{1}_{\omega}(\mrI)$, and 
$Y_{N} \, := \, \mbox{span}  \{ \what{G}_{j}^{(\beta \, , \, \alpha - \beta)} \}_{j = 0}^{N}
\subset H^{\alpha-1}_{\omega^*}(\mrI)$. Then a Petrov-Galerkin spectral approximation to 
\eqref{wform1} is: 
Given $f \in H^{-(\alpha-1)}_{\omega^{*}}(\mrI)$, determine $\phi_{N} \in X_{N}$
such that $u_{N}(x) \ = \ \omega(x) \, \phi_{N}(x)$ satisfies
\be
B(\phi_{N} \, , \, \psi_{N})
\ = \ \langle f \, , \, \psi_{N} \rangle_{\omega^{*}} \, , \ \ \forall \, \psi_{N} \in Y_{N} \, .
\label{appx1}
\ee
We first give the well posedness of \eqref{appx1}.

\begin{theorem} \label{exds1}
%
For $f \in H^{-(\alpha - 1)}_{\omega^{*}}(\mrI)$, $0 < k_{0} \le k(x)  \in L^{\infty}(\mrI)$, and $b \in W^{1 , \infty}_{w}(\mrI)$, 
$c \in L^{\infty}(\mrI)$ sufficiently small (with respect to $k_{0}$),
there exists a unique  $u_{N}(x) \ = \ \omega \, \phi_{N}(x)$ 
satisfying \eqref{appx1} with $\phi_{N}  \in X_{N}$,
and  $\| \phi_{N} \|_{H^{1}_{\omega}(\mrI)} \, \le \, \frac{1}{C_{3}} \| f \|_{H^{-(\alpha - 1)}_{\omega^{*}}(\mrI)}$,
for some positive constant $C_3$.
\end{theorem}
\textbf{Proof}: 
Similar to the proof of Lemma \ref{lmaBpos}, given $\phi_N\in X_N$ we can construct a $\psi_N\in Y_N$
such that there exists $C_{3} > 0$ such that
\begin{align}
\sup_{0 \ne \psi_{N} \in Y_{N}} 
\frac{ | B(\phi_{N} , \psi_{N}) |}{ \| \psi_{N} \|_{H^{\alpha-1}_{\omega^{*}}(\mrI)}} \ 
&\ge \ C_{3} \, \| \phi_{N} \|_{H^{1}_{\omega}(\mrI)} \, , \ \ \forall \, \phi_{N} \in X_{N} \, . \label{ddrtgf1}  
\end{align}
The proof then follows in an analogous manner to that of Theorem \ref{thmexun1}.
\mbox{ } \hfill \qed

Next we establish the error bound of $\phi-\phi_N$ in both the 
weighted $L^2$ and energy norms in the following theorem.
\begin{lemma} \label{errbd1}
Let $f \in H^{s}_{\omega^{*}}(\mrI)$ for $s \ge -(\alpha-1)$, 
$0 < k_{0} \le k(x)  \in W_{w}^{\ceil{s + \alpha - 1} , \infty}(\mrI)$, 
$b(x) \in W^{\max\{1 , \, \ceil{ \, \wtilde{s} \,} \}, \infty}_{w}(\mrI)$
and $c(x) \in 
W^{\ceil{ \, \wtilde{s} \,} , \infty}_{w}(\mrI)$ sufficiently small (with respect to $k_{0}$).
Then there exists $C > 0$ such that 
\begin{align}
 \| \phi \, - \, \phi_{N} \|_{H^{1}_{\omega}(\mrI)} &\ \le \ C \, N^{- (\wtilde{s} \, + \, \alpha-1)} \, 
 \| \phi \|_{H^{\wtilde{s} \, + \, \alpha}_{\omega}(\mrI)} \,
  \le \, C \, N^{- (\wtilde{s} \, + \, \alpha-1)} \, \| f \|_{H^{\wtilde{s}}_{\omega^{*}}(\mrI)}   ,
\label{herrest}\\
 \| \phi \, - \, \phi_{N} \|_{L^{2}_{\omega}(\mrI)}& \ \le \ C \, N^{- (\wtilde{s} \, + \, \alpha)} \, 
 \| \phi \|_{H^{\wtilde{s} \, + \, \alpha}_{\omega}(\mrI)} \,
  \le \, C \, N^{- (\wtilde{s} \, + \, \alpha)} \, \| f \|_{H^{\wtilde{s}}_{\omega^{*}}(\mrI)} \,  .
\label{L2errest}
\end{align}
\end{lemma}
\textbf{Proof}: We apply \eqref{ddrtgf1} to obtain for $\zeta_{N} \in X_{N}$
\begin{align}
C_{3} \,\| \phi_{N} \, - \, \zeta_{N} \|_{H^{1}_{\omega}(\mrI)}  
&\le \  \sup_{\stackrel{\psi_{N} \in Y_{N}}{\psi_{N} \neq 0}} 
\frac{ | B(\phi_{N} \, - \, \zeta_{N} \,  , \, \psi_{N}) |}{ \| \psi_{N} \|_{H^{\alpha-1}_{\omega^{*}}(\mrI)}}
= \  \sup_{\stackrel{\psi_{N} \in Y_{N}}{\psi_{N} \neq 0}} 
\frac{ | B(\phi \, - \, \zeta_{N} \,  , \, \psi_{N}) |}{ \| \psi_{N} \|_{H^{\alpha-1}_{\omega^{*}}(\mrI)}}
  \nonumber \\
&\le \sup_{\stackrel{\psi_{N} \in Y_{N}}{\psi_{N} \neq 0}} 
\frac{ C_{1} \, \| \phi \, - \, \zeta_{N}\|_{H^{1}_{\omega}(\mrI)} \, 
 \| \psi_{N} \|_{H^{\alpha-1}_{\omega^{*}}(\mrI)}}{ \| \psi_{N} \|_{H^{\alpha-1}_{\omega^{*}}(\mrI)}}  
\ = \  C_{1} \, \| \phi \, - \, \zeta_{N}\|_{H^{1}_{\omega}(\mrI)} .
\end{align}
Then, with $\zeta_{N} \, = \, P_{N} \phi \in X_{N}$, using Lemma \ref{lem:Approx}, 
\begin{align*}
  \| \phi \, - \, \phi_{N} \|_{H^{1}_{\omega}(\mrI)}
  &\le \ \| \phi \, - \, \zeta_{N} \|_{H^{1}_{\omega}(\mrI)} \ + \ 
    \| \zeta_{N} \, - \, \phi_{N} \|_{H^{1}_{\omega}(\mrI)}
\ \le \ C \,  \| \phi \, - \, \zeta_{N} \|_{H^{1}_{\omega}(\mrI)} \le \ C \, N^{- (\wtilde{s} \, + \, \alpha-1)} \, 
 \| \phi \|_{H^{\wtilde{s} \, + \, \alpha}_{\omega}(\mrI)}   \\
&\le \ C \, N^{- (\wtilde{s} \, + \, \alpha-1)} \, \| f \|_{H^{\wtilde{s}}_{\omega^{*}}(\mrI)} \, 
\   \mbox{ (using Corollary \ref{corregk1})}.
\end{align*}

To obtain the estimate for  $\|\phi \, - \, \phi_{N} \|_{L^{2}_{\omega}(\mrI)}$ we use
an  Aubin-Nitsche type argument. Introduce the following adjoint problem. 
Determine $\psi(x) \in  H^{\alpha - 1}_{\omega^{*}}(\mrI)$ satisfying
\[
\grave{\mcL}_{1-r}^{\alpha} (\omega^{*} \, \psi)(x) \  -  \ b(x) \, D (\omega^{*} \, \psi)(x) \ 
+ \ \big(c(x) \, - \, D b(x)) ( \omega^{*} \, \psi)(x) \ = \ (\phi \, - \, \phi_{N})(x) \, , \ x \in \mrI \, . 
\]
 As $\phi \, - \, \phi_{N} \in L^{2}_{\omega}(\mrI) \, = \, H^{0}_{\omega}(\mrI)$, from \eqref{kjg1} 
 we have
\be
  \| \psi \|_{H^{\alpha}_{\omega^{*}}(\mrI)} \ \le \ C \, \| \phi \, - \, \phi_{N} \|_{L^{2}_{\omega}(\mrI)} \, .
\label{ghre1}
\ee
Then, with  $\zeta_N \, = \, P_{N} \psi \, \in X_N$,
\begin{align*}
 \| \phi \, - \, \phi_{N} \|_{L^{2}_{\omega}} 
 &= \ \big( (\phi \, - \, \phi_{N}) \, , \, (\phi \, - \, \phi_{N}) \big)_{\omega} 
 \\
 &= \big( \phi \, - \, \phi_{N} \, , \, \grave{\mcL}^{\alpha}_{(1 - r)} \omega^{*} \, \psi \  -  \ b \, D \omega^{*} \, \psi 
 \ + \ \big(c \, - \, D b) \, \omega^{*} \, \psi  \big)_{\omega}  \\
& = \ B((\phi \, - \, \phi_{N}) \, , \, \psi)  = \ B((\phi \, - \, \phi_{N}) \, , \, \psi \, - \, \zeta_{N})  
\ \mbox{ (using Galerkin orthogonality)} \\
&\le \  C_{1} \, \| \phi \, - \, \phi_{N} \|_{H^{1}_{\omega}} \, \| \psi \, - \, \eta_{N} \|_{H^{\alpha-1}_{\omega^{*}}}  \\
&\le \ C \, N^{- (\wtilde{s} \, + \, \alpha-1)} \, 
 \| \phi \|_{H^{\wtilde{s} \, + \, \alpha}_{\omega}} \, N^{- 1} \, \| \psi \|_{H^{\alpha}_{\omega^{*}}}  
 \ \mbox{ (using \eqref{herrest} and Lemma \ref{lem:Approx})}  \\
&\le \ C \, N^{- (\wtilde{s} \, + \, \alpha)} \, 
 \| \phi \|_{H^{\wtilde{s} \, + \, \alpha}_{\omega}} \,  \| \phi \, - \, \phi_{N} \|_{L^{2}_{\omega}}  
  \ \mbox{ (using \eqref{ghre1})} \, ,
\end{align*}
from which \eqref{L2errest} follows.
\mbox{ } \hfill \qed

The weighted $L^2$ and energy estimates of $u - u_{N}$ immediately follow from 
Lemma \ref{errbd1} and Lemma \ref{lma4Her}.
\begin{corollary}  \label{uerrL2}
Let $f \in H^{s}_{\omega^{*}}(\mrI)$ for $s \ge -(\alpha-1)$, 
$0 < k_{0} \le k(x)  \in W_{w}^{\ceil{s + \alpha - 1} , \infty}(\mrI)$, 
$b(x) \in W^{\max\{1 , \, \ceil{ \, \wtilde{s} \,} \}, \infty}_{w}(\mrI)$
and $c(x) \in 
W^{\ceil{ \, \wtilde{s} \,} , \infty}_{w}(\mrI)$ sufficiently small (with respect to $k_{0}$).
Then there exists $C > 0$ such that 
\begin{align}
 \| u \, - \, u_{N} \|_{H^{1}_{\omega^{-1}}(\mrI)} \ 
&\le \  C \, N^{- (\wtilde{s} \, + \, \alpha-1)} \, \| f \|_{H^{\wtilde{s}}_{\omega^{*}}(\mrI)}  \, , \label{uHerrest}  \\
 \| u \, - \, u_{N} \|_{L^{2}_{\omega^{-1}}(\mrI)} \
&\le \  C \, N^{- (\wtilde{s} \, + \, \alpha)} \, \| f \|_{H^{\wtilde{s}}_{\omega^{*}}(\mrI)}  \, . \label{uL2errest}  
\end{align}
\end{corollary}

\subsection{Approximation of \eqref{wform1*}}
\label{sec_appx1*}
The main difference between the problems involving the different fractional diffusion operators,
$\acute{\mcL}_{r}^{\alpha}\cdot$ and  $\grave{\mcL}_{r}^{\alpha}\cdot$ is the functional spaces used.

Let
$X_{N} \, := \, \mbox{span}  \{ \what{G}_{j}^{(\alpha - \beta \, , \, \beta)} \}_{j = 0}^{N} 
\subset H^{\alpha - 1}_{\omega}(\mrI)$, and 
$Y_{N} \, := \, \mbox{span}  \{ \what{G}_{j}^{(\beta \, , \, \alpha - \beta)} \}_{j = 0}^{N}
\subset H^{1}_{\omega^*}(\mrI)$. Then a Petrov-Galerkin spectral approximation to 
\eqref{wform1} is: 
Given $f \in H^{-1}_{\omega^{*}}(\mrI)$, determine $\phi_{N} \in X_{N}$
such that $u_{N}(x) \ = \ \omega(x) \, \phi_{N}(x)$ satisfies
\be
B(\phi_{N} \, , \, \psi_{N})
\ = \ \langle f \, , \, \psi_{N} \rangle_{\omega^{*}} \, , \ \ \forall \, \psi_{N} \in Y_{N} \, .
\label{appx1*}
\ee

Analogous to Theorem \ref{exds1}, Lemma \ref{errbd1} and Corollary \ref{uerrL2} we have the following.
\begin{corollary}  \label{uerrL2*}
For $f \in H^{s}_{\omega^{*}}(\mrI)$, $-1 \le s$, 
$0 < k_{0} \le k(x)  
\in W_{w}^{\ceil{s + 1} , \infty}(\mrI)$, 
$b(x) \in W^{\max\{1 , \, \ceil{ \, \wtilde{s} \,} \}, \infty}_{w}(\mrI)$
and $c(x) \in 
W^{\ceil{ \, \wtilde{s} \,} , \infty}_{w}(\mrI)$
sufficiently small (with respect to $k_{0}$),
there exists a unique  $u_{N}(x) \ = \ \omega \, \phi_{N}(x)$ 
satisfying \eqref{appx1*} with $\phi_{N}  \in X_{N}$,
and  $\| \phi_{N} \|_{H^{\alpha - 1}_{\omega}(\mrI)} \, \le \, \frac{1}{C_{3}} \| f \|_{H^{-1}_{\omega^{*}}(\mrI)}$,
for some positive constant $C_3$. In addition there exists $C > 0$ such that
\begin{align}
 \| \phi \, - \, \phi_{N} \|_{H^{\alpha - 1}_{\omega}(\mrI)} &\ \le \ C \, N^{- (\wtilde{s} \, + \, 1)} \, 
 \| \phi \|_{H^{\wtilde{s} \, + \, \alpha}_{\omega}(\mrI)} \,
  \le \, C \, N^{- (\wtilde{s} \, + \, 1)} \, \| f \|_{H^{\wtilde{s}}_{\omega^{*}}(\mrI)}   ,
\label{herrest*}\\
 \| \phi \, - \, \phi_{N} \|_{L^{2}_{\omega}(\mrI)}& \ \le \ C \, N^{- (\wtilde{s} \, + \, \alpha)} \, 
 \| \phi \|_{H^{\wtilde{s} \, + \, \alpha}_{\omega}(\mrI)} \,
  \le \, C \, N^{- (\wtilde{s} \, + \, \alpha)} \, \| f \|_{H^{\wtilde{s}}_{\omega^{*}}(\mrI)} \,  .
\label{L2errest*}  \\
 \| u \, - \, u_{N} \|_{H^{\alpha - 1}_{\omega^{-1}}(\mrI)}  
&\ \le \  C \, N^{- (\wtilde{s} \, + \, 1)} \, \| f \|_{H^{\wtilde{s}}_{\omega^{*}}(\mrI)}  \, , \label{uHerrest*}  \\
 \| u \, - \, u_{N} \|_{L^{2}_{\omega^{-1}}(\mrI)} 
&\ \le \  C \, N^{- (\wtilde{s} \, + \, \alpha)} \, \| f \|_{H^{\wtilde{s}}_{\omega^{*}}(\mrI)}  \, . \label{uL2errest*}  
\end{align}
\end{corollary}

\setcounter{equation}{0}
\setcounter{figure}{0}
\setcounter{table}{0}
\setcounter{theorem}{0}
\setcounter{lemma}{0}
\setcounter{corollary}{0}
\setcounter{definition}{0}
%
\section{Model Discussion and Numerical Experiments}
\label{sec_num}

In this section we carry out numerical experiments to demonstrate the accuracy and performance of the scheme and to study the behavior of the two models.

\subsection{Accuracy test}
In this subsection we test the accuracy of the scheme (\ref{appx1}) under different parameters 
\begin{itemize}
\item[(a)] $\alpha = 1.30$, $r = 0.50$, $k(x)=1+2x$;
\item[(b)] $\alpha = 1.60$, $r = 0.40$, $k(x)=1-0.3\sin x$.
\end{itemize}
  The other parameters are selected as those in \cite[Experiment 2]{zhe211}: $b(x) = e^x$, $c(x) =  5+\sin(x)$, and $f(x) = 1$.
Then we could follow the results in Corollary \ref{uerrL2} to predict the convergence rates for the errors under both the weighted $L^2$ and $H^1$ norms. The numerical
convergence rates are presented in Tables \ref{table1}-\ref{table2}, which are in good agreement
with the predicted rates.

\begin{table}[h!]
	\setlength{\abovecaptionskip}{0pt}
	\centering
	\caption{Accuracy test for case (a)}	\label{table1}
	\vspace{0.5em}	
	\begin{tabular}{ccccc}
		\hline
		$N$&$\|u-u_N\|_{L^2_{\omega^{-1}}}$ &$\kappa$& $\|u-u_N\|_{H^{1}_{\omega^{-1}}}$ &$\kappa$ \\
		\cline{1-5}
8&	6.50E-03&		&6.81E-02	&\\
10&	4.19E-03&	1.96 &	5.31E-02&	1.12\\ 
12& 2.91E-03&	2.01 &	4.28E-02&	1.18 \\
14&	2.12E-03&	2.04 &	3.54E-02&	1.24 \\
16&	1.62E-03&	2.05 &	2.97E-02&	1.30 \\
		\hline	
		Pred.&    &2.25&     &  1.25 \\
		\hline	
	\end{tabular}
\end{table}

\begin{table}[h!]
	\setlength{\abovecaptionskip}{0pt}
	\centering
	\caption{Accuracy test for case (b)}	\label{table2}
	\vspace{0.5em}	
	\begin{tabular}{ccccc}
		\hline
		$N$&$\|u-u_N\|_{L^2_{\omega^{-1}}}$ &$\kappa$& $\|u-u_N\|_{H^{1}_{\omega^{-1}}}$ &$\kappa$ \\
		\cline{1-5}
8&	3.37E-03&		&3.61E-02	&\\
10&	1.70E-03&	3.06 &	2.21E-02&	2.19\\ 
12&	9.70E-04&	3.08 &	1.49E-02&	2.19 \\
14&	6.09E-04&	3.03 &	1.07E-02&	2.15 \\
16&	4.10E-04&	2.96 &	8.04E-03&	2.11 \\
		\hline	
		Pred.&    &2.95&     &  1.95 \\
		\hline	
	\end{tabular}
\end{table}

\subsection{Model comparison}

In this subsection we carry out numerical experiments to study the behavior of the two models.

\begin{figure}[h!]
	\setlength{\abovecaptionskip}{0pt}
	\centering
	\includegraphics[width=3.2in,height=3.2in]{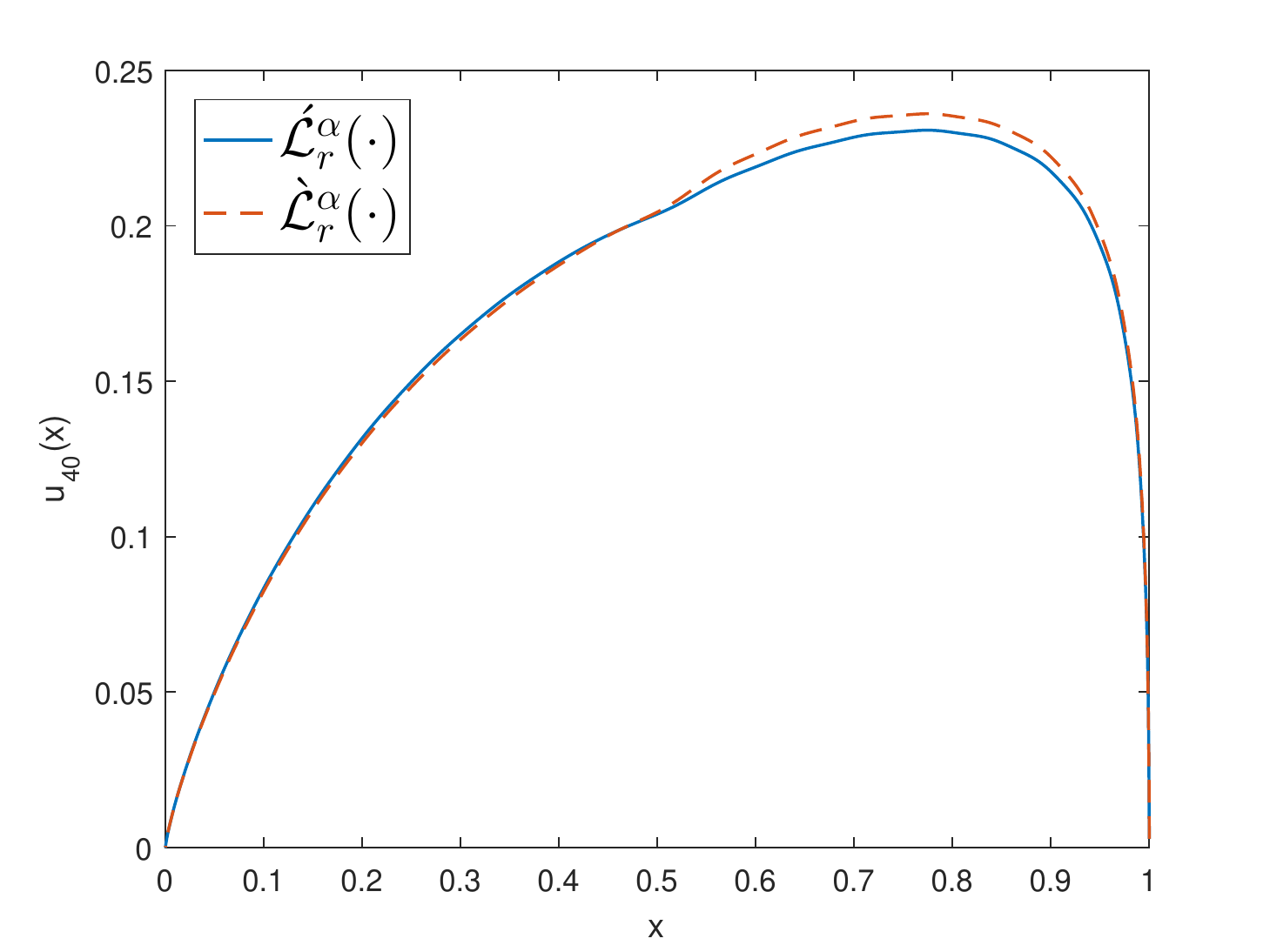}
	\includegraphics[width=3.2in,height=3.2in]{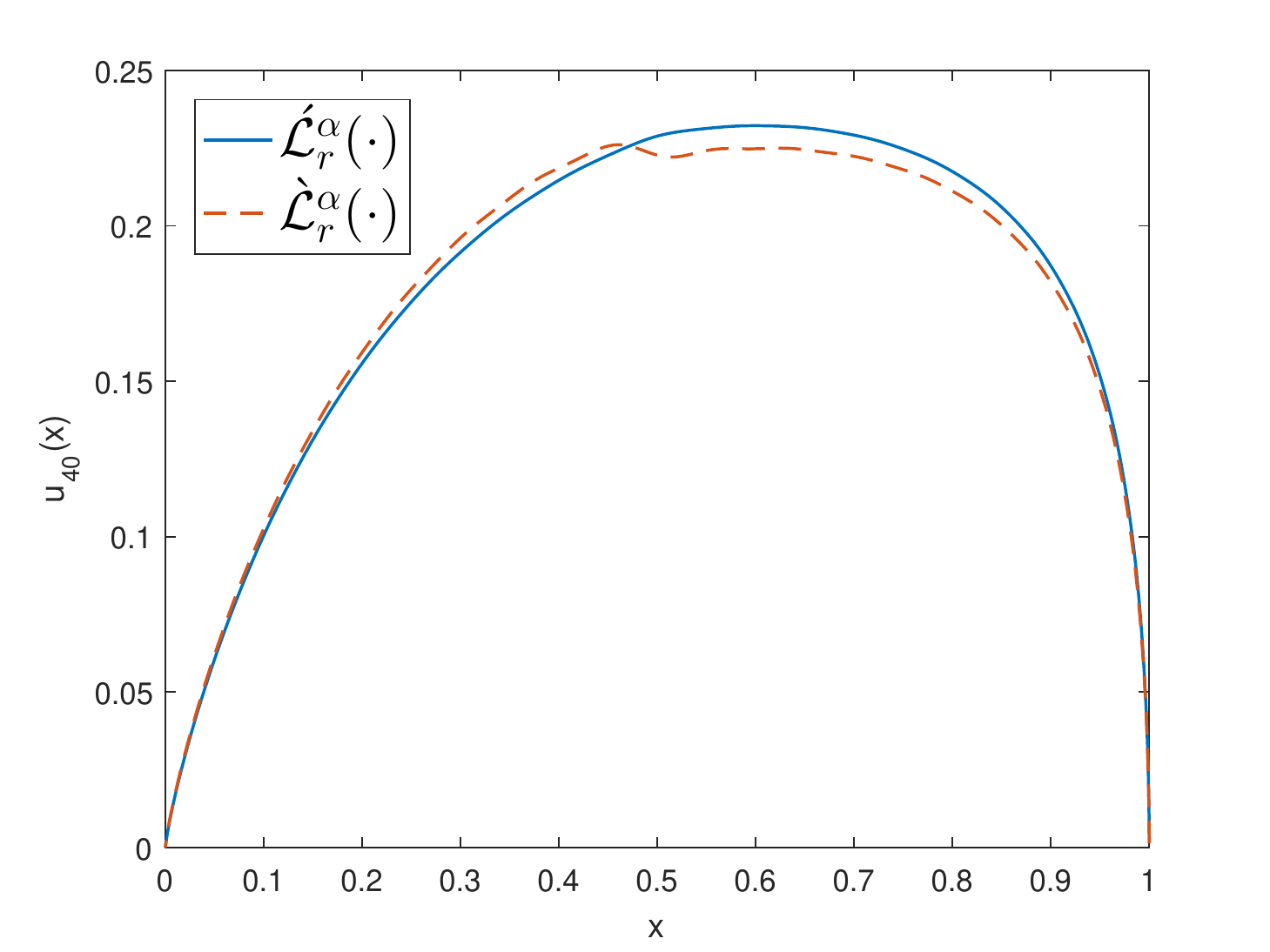}
	\caption{Plots of the reference solution $u_{40}(x)$ under (left) $k_1(x)$ and (right) $k_2(x)$.}
	\label{plot1}
\end{figure}

We compare the solutions for model (\ref{DefProb2})-(\ref{DefBC2}) involving both $\acute{\mcL}_{r}^{\alpha}(\cdot)$ and $\grave{\mcL}_{r}^{\alpha}(\cdot)$. Let $\alpha=1.4$, $r=0.4$, $b(x)$, $c(x)$ and $f(x)$ are given as before, and we take the practical diffusivity coefficient $k(x)$ as follows
\begin{equation}\label{defmcL2:e0}
  k_1(x) \ = \ \left\{ \begin{array}{rl}
                      2, & \mbox{ for } 0 \leq x < 1/2; \\
                      1, & \mbox{ for } 1/2 \leq x \leq 1 \, ,
                      \end{array} \right.
                     ~~
  k_2(x) \ = \ \left\{ \begin{array}{rl}
                      1, & \mbox{ for } 0 \leq x < 1/2; \\
                      2, & \mbox{ for } 1/2 \leq x \leq 1 \, ,
                      \end{array} \right.
\end{equation}                      
and present the numerical results in Figure \ref{plot1}. 

To better understand these results we note that the operator $\grave{\mcL}_{r}^{\alpha}u(x)$ in \eqref{defmcL2} can be expressed as a local mass balance with a two-sided fractional Fick's law
\begin{equation}\label{defmcL2:e1}
\grave{\mcL}_{r}^{\alpha}u(x) = D \grave{\mathcal F}^\alpha(x), ~~\mbox{with}~\grave{\mathcal F}^\alpha(x) :=  - k(x) \big( r D^{-(2 - \alpha)} \ + \ (1 - r) D^{-(2 - \alpha)*}  \big) \, D u(x). 
\end{equation}
Hence, the model \eqref{DefProb2} and \eqref{defmcL2:e1} is a nonlocal extension of the classical local mass balance with a local Fick's law. Like its second-order analogue, the interface conditions for the strong solutions of problem \eqref{DefProb2}, \eqref{DefBC2}, and \eqref{defmcL2} now take the form
\begin{equation}\label{defmcL2:e2}
u\big((1/2)_- \big) = u\big((1/2)_+ \big), \quad \grave{\mathcal F}^\alpha\big((1/2)_- \big) = \grave{\mathcal F}^\alpha\big((1/2)_+\big).
\end{equation}
The combination of the continuity of the fractional diffusive flux in \eqref{defmcL2:e2} and the jump discontinuity of the material diffusivity coefficients $k_1(x)$ and $k_2(x)$ in \eqref{defmcL2:e0} across the interface, and the fact that $\grave{\mathcal F}^\alpha(x)$ depends only on $k(x)$ enforce that the two-sided fractional derivative of $u$ has a jump discontinuity across the interface $x=1/2$, which is proportional to the jump size of the diffusivity coefficient $k_1$ and $k_2$ across the interface, although $u$ is continuous across the interface. The continuity of $u$ and the jumps of the fractional drivative of $u$ across the interface can be observed clearly in Figure \ref{plot1}.

On the other hand,  the operator $\acute{\mcL}_{r}^{\alpha}u(x)$ in \eqref{defmcL} can be expressed as a local mass balance with a nonlocal Fick's law, in which the flux at a point is a weighted average of local fluxes throughout the 
domain
\begin{equation}\label{defmcL:e1}
\acute{\mcL}_{r}^{\alpha}u(x) = D \acute{\mathcal F}^\alpha(x), ~~\mbox{with}~\acute{\mathcal F}^\alpha(x) :=  \big( r D^{-(2 - \alpha)} \ + \ (1 - r) D^{-(2 - \alpha)*}  \big) \big (-k(x) \, D u(x)\big). 
\end{equation}
In this case the interface conditions 
\begin{equation}\label{defmcL:e2}
u\big((1/2)_- \big) = u\big((1/2)_+ \big), \quad \acute{\mathcal F}^\alpha\big((1/2)_- \big) = \acute{\mathcal F}^\alpha\big((1/2)_+\big)
\end{equation}
assume the same form as \eqref{defmcL2:e2} but with a different flux function $\acute{\mathcal F}^\alpha$ defined in \eqref{defmcL:e1}. In contrast to $\grave{\mathcal F}^\alpha$, at each $x \in (0,1)$ $\acute{\mathcal F}^\alpha(x)$ depends on $k(y)$ for all the $y \in (0,1)$. In particular, the continuity of $u$ and the flux $\acute{\mathcal F}^\alpha$ across the interface specified in \eqref{defmcL:e2} imply
\begin{equation}\label{defmcL:e3}\begin{array}{rl}
\acute{\mathcal F}^\alpha\big((1/2)_- \big) & \displaystyle = - \bigg \{ \frac{r}{\Gamma(2 - \alpha)}   \int_{0}^{\frac{1}{2}} \frac{k(s)u'(s)}{(\frac{1}{2} - s)^{\alpha - 1}} ds + \frac{1-r}{\Gamma(2 - \alpha)} \int_{\frac{1}{2}}^{1} \frac{k(s)u'(s)}{(s - \frac{1}{2})^{\alpha - 1}} ds \bigg\} \\
& \displaystyle = \acute{\mathcal F}^\alpha\big((1/2)_+\big)
\end{array}\end{equation}
which holds naturally without enforcing a jump discontinuity of the derivative of $u$ across the interface, as shown in Figure \ref{plot1}.

\section{Conclusions}
In this paper we prove the well-posedness of the the variable coefficient two-sided fractional diffusion, advection, reaction equations on a bounded interval, which remains untreated in the literature. Additionally we establish the regularity of
the solution in terms of the regularity of the coefficient functions and the right hand side function. A Petrov-Galerkin
scheme for the approximation of the solution is proposed and analyzed.
A main contribution of this work lies in designing appropriate test and trial functions to prove the inf-sup condition of the variable coefficient fractional diffusion, advection, reaction operators in suitable function spaces, upon which the subsequent mathematical and numerical analysis rely. Numerical experiments are presented to substantiate the theoretical findings and to compare the behaviors of different models.


\section*{Declarations}
 This work was partially supported by the ARO MURI Grant W911NF-15-1-0562, by the National Science Foundation under
Grant DMS-2012291, by the China Postdoctoral Science Foundation under Grants 2021TQ0017 and 2021M700244, by the National Natural Science Foundation of China under Grant 12071262, and by International Postdoctoral Exchange
Fellowship Program (Talent-Introduction Program) YJ20210019.
  
 All data generated or analyzed during this study are included in this published article.


\end{document}